\documentclass[12pt]{article}
\begin{document}

\title{Global wellposedness and scattering for 3D Schr\"odinger equations with harmonic potential and radial data }
\author{Zhang Xiaoyi$^{1,2}$
\thanks{The author is supported in part by N.S.F. grant
No.10426029 and the Morningside Center.}\\
$^1$Institute of Mathematics , Chinese Academy of Sciences,\\
 $^2$Beijing Institute of Applied Physics and Computational Mathematics\\
    P.O.Box 8009-11, Beijing 100088, China\\
 Email:xiaoyizhangoo@yahoo.com.cn
}
\date{}
\maketitle

\abstract{ In this paper,we show that spherical bounded energy
solution of the defocusing 3D energy critical Schr\"odinger
equation with harmonic potential, $(i\partial_t + \frac
{\Delta}2+\frac {|x|^2}2)u=|u|^4u$, exits globally and scatters to
 free solution in the space $\Sigma=H^1\bigcap\mathcal F H^1$.
We preclude the concentration of energy in finite time by
combining the energy decay estimates and the ideas in the
paper\cite{Bourgain}\cite{CKSTT}\cite{Tao}.}



\section{Introduction}
In this paper, we consider the Cauchy problem of defocusing energy
critical equation with harmonic potential,
\begin{eqnarray}
&& (i\partial_t +\frac {\Delta}2)u=-\frac {|x|^2}2
u+|u|^4u,\qquad (x,t)\in \mathbf R^3\times \mathbf R ,\label{1.1}\\
&&u(x,0)=u_0(x),\qquad x\in \mathbf R^3 , \label{1.2}
\end{eqnarray}
where $u(t,x)$ is a complex function on $\mathbf R^3\times \mathbf
R$, $u_0(x)$ is a complex function on $\mathbf R^3$ satisfying
$$
u_0(x)\in \Sigma =\{v;
\|v\|_{\Sigma}=\|v\|_{H^1}+\|xv\|_2<\infty\}.
$$
We will be interested in the global existence and long time
behavior of the solution.

Schr\"odinger equation without potential has been extensively
studied and we are mainly interested in the equation with power-like nonlinearity:
\begin{equation}
(i\partial_t+\frac {\Delta}2)v=\mu |v|^pv, \qquad
\mu\neq 0, \label{2}
\end{equation}
where $\mu>0$ and $\mu<0$ correspond to the defocusing case and
the focusing case, respectively. We concern such Cauchy problem
as: there exists $t_0\in \mathbf R$ such that the solution $v$ is
in $H^1$ at this point. One thing that plays important role in the
study of the Cauchy problem is the conservation of energy:
$$
E(t)=\frac12\|\nabla v(t)\|_2+\frac {\mu}{p+2}\|v(t)\|_{p+2}^{p+2}=const.
$$
Roughly speaking, (\ref{2}) has local $H^1$ solution when $p$ is
smaller than and equals to a certain exponent $p_c$ ($p_c=4$ as
$n=3$) which is energy critical in the sense that the natural
scale invariance
$$
v(x,t)\longrightarrow \lambda ^{-\frac 2p} v(\lambda ^{-1}x,\lambda^{-2}t)
$$
of the equation leaves the $\dot H^1$ norm invariant. In the
supercritical case $p>p_c$, (\ref{2}) is locally illposed in the sense that the solution $v(x,t)$ does not continuously
 depend on the initial data $v_0(x)$ in $H^1$ space. To get more details, one refers to
 see \cite{Cazenave1}, \cite{CW}, \cite{GV1},
 \cite{CCT}.
In the defocusing and subcritical case ,the global existence is a
direct consequence of the energy conservation. In the focusing
case $\mu<0$, blow up in finite time may appear, since the
influence of kinetic energy $ \|\nabla v\|_2^2$ may not always
surpass the influence of the potential energy $
\|v\|_{p+2}^{p+2}$. See \cite{Glassey}, \cite{Cazenave1} for
instance. The case of energy critical becomes rather difficult
because the pure energy conservation is not enough to ensure the
energy solution exist globally. In other words, the energy of the
solution may concentrate somewhere, so that the solution may
possibly blow up in finite time.

The first important work in this field is due to
Bourgain\cite{Bourgain} and Grillarkis \cite{Grillakis2}. They
pointed out that the solution will concentrate somewhere in
$\mathbf R^d$ unless the solution does not exist globally. To
preclude this phenomenon, they use an apriori estimate which is
called Morawetz estimate,
\begin{equation}
\int_I \int _{A|I|^{\frac 12}} \frac {|v|^{p_c+2}(x,t)}{|x|} dxdt
\le CA E(v) |I|^{\frac12},\quad A\ge 1.\label{3}
\end{equation}
(\ref{3}) is useful for preventing the concentration of $v(t,x)$
at origin $x=0$. This is especially helpful when the solution is
radially symmetric since in this case, it can be easily shown by
bounded energy that $v$ will not concentrate at any other location
than the origin. Their results are restricted to $d=3,4$. Recently
in\cite{Tao}, T.Tao proved the energy critical Cauchy problem with
radial data is globally wellposed and scatters to free solution in
all dimensions by a different method.

To remove the radial assumption is not an easy thing. Recently in \cite{CKSTT}, J. Colliander and others
proved this is feasible in the 3-dimensional
case. Their proof involves many technical analysis especially in both physical space and frequency space. For details see
\cite{CKSTT}.

Schr\"odinger equation with harmonic potential and power-like
nonlinearity can be written in the form,
 \begin{equation}
 (i\partial_t+\frac{\Delta}2)u=\frac{\omega |x|^2}2u+\mu |u|^pu.
 \label{4}
 \end{equation}
The Cauchy problem of it becomes much more complicated because of
many-sided reasons. One of them lies in getting the estimate of
the linear operator, another one may be that general structural
conditions such as scaling invariance and spatially translation
invariance that hold for the equation (\ref{2}) will not hold for
(\ref{4}) any more. Even so, when we study Cauchy problem, we
still define the energy critical exponent by omitting the
potential just like we do for the equation (\ref{2}). In both
subcritical and critical case, it's natural to seek for the
solution with suitable decay in space, ie,
\begin{equation}
u\in C_t(\Sigma).
\end{equation}
Indeed, recently in \cite{Carles1}, \cite{Carles2}, R. Carles
systematically studied the Cauchy problem of (\ref{4}). He found
that when the nonlinearity is subcritical, defocusing and when a
confining potential (ie.$\omega=-1$) is involved, the solution
will be global and a scattering theory is available. Furthermore,
when the nonlinearity is focusing and subcritical, a sufficient
strong harmonic potential will prevent blow up in finite time.
There are also some results involving a general potential $V(x)$
which is a quadratic function.

There are two things that play important role in his proof. The
first ingredient involved the Strichartz estimate for the linear
operator $i\partial_t+\frac{\Delta}2+\frac{\omega|x|^2}2$. This
may follow from Mehler's formula at some special occasions. The next thing is that
there exist two Galilean  operators $J(t)$ and $H(t)$ which can commute with the
linear operator and can be viewed as a substitute of
$\nabla$ and $x$ in the nonpotential case.

The question remains open that what will happen when the
nonlinearity is energy critical, that is, $p=p_c=\frac 4{d-2}$,
$d\ge 3$. In this paper, we restrict our attention to the case
$d=3$, $\omega=-1$, and $\mu=1$. Enlightened by the work of
Bourgain, Carles, Tao and others, we are expected to prove the
global existence and scattering theory for the cauchy problem
(\ref{1.1})-(\ref{1.2}). In what follows, we sketch the proof.

The first thing we do is to find the local solution and small
global solution of (\ref{1.1}) and (\ref{1.2}). They are available
thanks to the Strichartz estimate for the linear operator. One
interesting thing in the small solution theory is that: in order
to get global solution, it suffices to require $\|\nabla u_0\|_2$
to be small enough. This result is fundamental in our proof.
Another thing we must pay attention is that the maximal time
interval of the local solution depends on the profile of the
initial data, not depends on the $\Sigma$ norm of $u_0$ only. This
is the main reason why the critical problem becomes much more
complicated than the subcritical case.

Next, we show that: in order to extend the local solution to a global one and prove
 the scattering theory to the global solution, it suffices to
 prove an apriori space time bound for the solution. Being same
 with \cite{Bourgain}, this bound is $\|u\|_{L^{10}(R;L^{10})}$. So, it's natural for us to
 see what the apriori bounds have been provided by the equation
 (\ref{1.1}). First of all, we have two conservative quantities:
 mass and energy,
\begin{eqnarray}
&& \|u(t_1)\|_2=\|u(t_2)\|_2,  \\
\label{5} && E(u(t_1))= E(u(t_2))=\frac 12\|\nabla
u(t_2)\|_2^2-\frac 12\|xu(t_2)\|_2^2+\frac13\|u(t_2)\|_6^6,
\label{6.0}
\end{eqnarray}
the energy (\ref{6.0}) is non-positive, so we split it into two
positive parts:
\begin{eqnarray}
\begin{array}{ll}
E_1(u(t))=\frac 12\|\nabla u(t)\|_2^2+\frac 13\|u(t)\|_6^6,\\
E_2(u(t))=\frac12 \|x u(t)\|_2^2,
\end{array}
\end{eqnarray}
and consider the Cauchy problem with $E_1(u(0))=E$, $E_2(u(0))=B$ for some fixed constant $E>0$, $B>0$.
 From now on, we call a finite energy solution $u(x,t)$ on a time interval
$I$ if it is such that $E_i(u(t))<\infty$, $\forall t\in I$,
$i=1,2$, and we may write the notations by omitting $u$ in some
occasions for the sake of simplicity. Although $E_1(t)$ and
$E_2(t)$ are nonnegative all the time, we are not clear about the
evolution of them. So we introduce another way that is provided by
R. Carles\cite{Carles1}, \cite{Carles2}:
\begin{eqnarray}
\begin{array}{ll}
\mathcal E_1(t)=\frac 12\|J(t)u(t)\|_2^2+\frac 13\cosh
^2t\|u(t)\|_6^6,\\
\mathcal E_2(t)=\frac 12\|H(t)u(t)\|_2+\frac
13\sinh^2t\|u(t)\|_6^6,
\end{array}
\end{eqnarray}
$\mathcal E_1(t)-\mathcal E_2(t)=E(t)$ and they coincide with
$E_1(t)$ and $E_2(t)$ only at $t=0$.
 The benefits of this decomposition is that we know that
 $\mathcal E_1(t)$, $\mathcal E_2(t)$ all decay in time. (See
 Section2 for details). Using this facts combing the Strichartz
 estimate, it's not difficult to get global solution in the
 subcritical case.

In the critical case, the decay estimates of $\mathcal E_1(t)$ and $\mathcal E_2(t)$ are not sufficient to prevent blow
up in finite time. However, it's helpful in the sense that it provides strong decay of the potential energy $\|u(t)\|_6^6$.
From this and some elementary analysis, we can  fix a time $T$ only dependent with $E$ such that on $(-\infty, -T)\cup
(T,\infty)$, $\|u\|_{10}$ has good control. Thus, we are left to do estimates on a finite time interval $[-T(E),T(E)]$.

By the decay estimates of $\mathcal E_1(t)$ and $\mathcal E_2(t)$ and the relation between $E_i(t)$ and
$\mathcal E_i(t)$, we know that on $[-T(E),T(E)]$, $E_1(t)$ and $E_2(t)$
are bounded uniformly by constants $\Lambda_1(E,B)$, $\Lambda_2(E,B)$ respectively. Now, we fix a small constant $\eta_1=
\eta_1(\Lambda_1,\Lambda_2+\Lambda_1)$ and divide $[-T(E),T(E)]$ into finite intervals with fixed length $\eta_1^4$. If we can prove
that on each subinterval, $u$ has finite $L^{10}$ estimate bounded only by $C(\Lambda_1,\Lambda_2)$, we can sum these
intervals together and give the final result.

Now, let's clarify again what is left to do. Let $t_0\in
\mathbf{R}$, and $v(t_0)\in \Sigma$ satisfy $E_1(v(t_0))\le
\Lambda_1$, and $E_2(v(t_0))\le \Lambda_2$, then we are required
to prove that there exists constant
$\eta_1(\Lambda_1,\Lambda_2+\Lambda_1)$ such that the Cauchy
problem with the prescribed data $v(t_0)$ at time $t_0$ is at
least solvable on $[t_0-\eta_1^4, t_0+\eta_1^4]$ and satisfies the
estimate:
$$\|v\|_{L^{10}([t_0-\eta_1^4,t_0+\eta_1^4];L^{10})}\le C(
\Lambda_1,\Lambda_2).$$

Thanks to the local solution theory, we need only to prove the
above estimate by apriorily  assuming that the solution has
existed on interval $[t_0-\eta_1^4,t_0+\eta_1^4]$. Here, we adopt
the ideas in \cite{Bourgain} to get this estimate.

 By time translation, we may assume $t_0=0$. Fix the small constant
$\eta_1$ such that it satisfy all the conditions that will appear
in the proof, we subdivide $[0,\eta_1^4]$ into $J_1$ intervals and
$[-\eta_1^4, 0]$ into $J_2$ intervals such that on each
subinterval, $v$ has $L^{10}$ norm comparable with $\eta_1$. We do
analysis forward in time and aim to estimate $J_1$ for simplicity.
By some technical computation and the radial assumption, we get a
sequence of bubbles located at the origin for a sequence of times
in each subinterval. If the volume of every bubble is sizeable by
the length of the corresponding time interval, then the solution
is solitonlike and $J_1$ can be estimated by using Morawetz
estimate. Otherwise, there is concentration for $E_1(v(t_*))$ for
some $t_*\in(0,\eta_1^4)$. Our main task is to estimate $J_1$ in
this case.

By removing the small bubble(because of the concentration), we get
a new function $w(t_*)$ for which $E_1(w(t_*))\le E_1
(v(t_*))-c\eta_1^3$. Here, we meet with a problem in making
comparison between $E_1(v(t_*))$ and $E_1(v(0))$ because as we
have mentioned before, we are not clear about the evolution of
$E_1(t)$ in time. However, thanks to the previous simplification
to the initial problem, and by using the small length condition,
we are able to roughly estimate the increment of $E_1(v(t))$ from
0 to $t_*$ by $C\eta_1^4$. Therefore, we get the final estimate
$E_1(v(t_*))\le \Lambda_1-c\eta_1^3$. This allows us to do
induction on the size of energy $E_1$. Another difficulty comes
from $E_2$ since there is no concentration property for it.
However, we can deal with this trouble by noticing that the
increment from $E_2(v(0))$ to $E_2(w(t_*))$ is also small, and any
finite increment during the iteration is permitted by the small
solution theory. This is the reason why we take
$\eta_1=\eta_1(\Lambda_1,\Lambda_1+\Lambda_2)$ instead of
$\eta_1(\Lambda_1, \Lambda_2)$. By considering all the factors
together, we may make an inductive assumption as follows:

 Let $v(t')\in \Sigma$ satisfy $E_1(v(t'))\le \Lambda_1-C\eta_1^4$
and $E_2(v(t'))\le \Lambda_2+C\eta_1^4$, then the Cauchy problem
of (\ref{1.1}) with prescribed data $v(t')$ at time $t'$ is at
least solvable on $[t'-\eta_1^4,t'+\eta_1^4]$, and there holds
that $\|v\|_{L^{10}([t'-\eta_1^4, t'+\eta_1^4];L^{10})}\le
C(\Lambda_1-C\eta_1^4, \Lambda_2+C\eta_1^4)$.

By using this assumption, $J_1$ and $J_2$ can be estimated by some technical arguments.

Finally, Let's explain why we do induction on the size of $E_1(t)$
and $E_2(t)$ and not on the size of $\mathcal E_1(t)$ and
$\mathcal E_2(t)$, since at first glance, the latter has good
decay, thus is hopeful to be viewed as a substitute of Hamiltonian
for Schr\"odinger equation without potential. Another reason
supporting the idea is that, one can get small solution once for
some $t\in R$, $\|J(t)u(t)\|_2^2$ is sufficiently small.
 However, we notice that, not liking the quantity $E_i(t)$, the quantity $\mathcal E_i(t)$ is not time-translation
 invariant, this will make essential trouble and is the key reason that one should not do
 induction on the size of $\mathcal E_i(t)$.

The remaining part of this paper is arranged as follows: In
Section2, we give some notations and some basic estimates. They
include: Littlewood-Paley decomposition, Galilean operator,
Strichartz estimates for the linear operator with potential, basic
properties of Galilean operator, etc. In the first part of
Section3, we give the local wellposedness and small solution
theory. The small solution theory claims that the Cauchy problem
of (\ref{1.1}) is global
 wellposed
and scatters to free solution if for some $t_0$, $E_1(t_0)$ is small enough. This is the fundamental theory which allows
us to do induction. In the second part, we use the decay estimate to simplify the large data problem to an aproiri estimate
on a finite time interval. Section4 is devoted to Morawetz estimate of the solution of (\ref{1.1}). In Section5, we use
Littlewood-Paley and paraproduct decomposition to prove the existence of a sequence of bubbles.
In Section6, we control $J_1$ and $J_2$ in the case
of solitonlike solution. In Section7, We control $J_1$ and $J_2$ if there is concentration by using the inductive
assumption and close the induction by a perturbation analysis in Section8.

\section{Notations and basic estimates}

\noindent Notations:

Let $\eta_1$, $\eta_2$, $\eta_3$ be small numbers satisfying
$0<\eta_3<\eta_2<\eta_1$ and to be defined in the proof,
$c(\eta_1)$, $c(\eta_2)$, $c(\eta_3)$ be small numbers satisfying
$0<c(\eta_3)<c(\eta_2)<c(\eta_1)\ll 1$; $C(\eta_1)$, $C(\eta_2)$,
$C(\eta_3)$ be large numbers such that $1\ll C(\eta_1)\ll
C(\eta_2)\ll C(\eta_3)$. $C$, $c$ are absolute numbers and may be
different from one line to another.

For any time interval $I$, we use $\|\cdot\|_{L^q(I;L^r)}$ to denote the Lebesgue norm, where $1\le q, r\le \infty$.

Next, we give the definition of Littlewood-Paley projection. Let
$\{\phi_j(\xi)\}_{j=-\infty}^{j=\infty}$  be a sequence of smooth
functions and each supported in an annuli $\{\xi;
2^{j-1}\le|\xi|\le 2^{j+1}\}$, furthermore, for any $\xi\neq0$,
$$
\sum_{j=-\infty}^{\infty}\phi_j(\xi)=1.
$$
For any $N=2^j$, we define Littlewood-Paley projection as follows:
\begin{eqnarray*}
&&P_N=P_{2^j}=\mathcal F^{-1}(\phi_j)*\cdot,\\
&&P_{\le N}=P_{\le 2^j}=\mathcal F^{-1}(\sum_{j'\le j}\phi_{j'})*\cdot,\\
&&P_{>N}=I-P_{\le N}.
\end{eqnarray*}
We list some basic properties of the projector which will be used
often :

$\cdot$ For any $1\le p\le\infty$, and $s\ge 0$, we have:
\begin{eqnarray*}
&&\||\nabla|^s P_N f\|_p\sim N^s\|P_N f\|_p,\\
&&\||\nabla|^s P_{\le N}\|_p\le C N^s\|P_{\le N}\|_p,
\end{eqnarray*}

$\cdot$ Beinstein estimate: For any $1\le q\le p\le\infty$, we have
\begin{eqnarray*}
\|P_Nf\|_p\le CN^{d(\frac 1q-\frac 1p)}\|P_Nf\|_q, \\
\|P_{\le N}f\|_p\le CN^{d(\frac 1q-\frac 1p)}\|P_{\le N}f\|_q.
\end{eqnarray*}

Let $u(t,x)$ be the solution of 3-d linear Schr\"odinger equation
with confining potential:
\begin{equation}
\begin{array}{ll}
(i\partial_t+\frac{\Delta}2)u=-\frac{|x|^2}2 u,\\
u|_{t=0}(x)=u_0(x),
\end{array}
\label{6}
\end{equation}
then it can be expressed by the Mehler's formula (see
\cite{FH65}),
\begin{eqnarray}
\displaystyle u(t,x)&=&U(t) u_0=e^{-\frac {it}2(-\Delta-|x|^2)}u_0 \nonumber \\
\displaystyle &=&e^{-\frac{i3\pi}4 sgn t} |\frac 1{2\pi \sinh
t}|^{\frac 32} \int_{\mathbf R^3} \displaystyle e^{\frac i{\sinh
t}(\frac {x^2+y^2}2 \cosh t -x\cdot y)} u_0(y) dy, \label{7}
\end{eqnarray}
one sees from the above that the kernel of $U(t)$ has the better
dispersive estimate than the kernal of Schr\"odinger operator
without potential. By using Mehler's formula (\ref{7}), and noting
that $U(\cdot)$ is unitary on $L^2$, one has the following decay
estimate
\begin{eqnarray}
&&\|U(t)u_0\|_{\infty}\le C |t|^{-\frac32}\|u_0\|_1,
\label{7.1} \\
&&\|U(t)u_0\|_p\le
C|t|^{\frac32(\frac1p-\frac1{p'})}\|u_0\|_{p'},\quad 2\le
p\le\infty .\label{7.2}
\end{eqnarray}
Using this decay estimate and by some standard arguments, one can
get Strichartz estimates for the operator $U(t)$.

{\bf {Definition 2.1}} A pair $(q,r)$ is admissible if $2\le r<6$ and $\displaystyle\frac 2q +\frac 3r=\frac 32$.

{\bf Lemma 2.2}\ Strichartz estimates for $U(t)$.

For any admissible pair $(q,r)$, there exists $C_r>0$ such that
$$
\|U(\cdot) \phi\|_{L^q(\mathbf R;L^r)}\le C_r\|\phi\|_2.
$$

For any admissible pairs $(q_1,r_1)$, $(q_2,r_2)$ and any time
interval $I$, there exists constant $C_{r_1,r_2}$,  such that
$$
\|\int_{I\cup \{s<t\}} U(t-s)F(s)ds\|_{L^{q_1}(I;L^{r_1})} \le C_{r_1r_2} \|F\|_{L^{q_2'}(I;L^{r_2'})}.
$$
We omit the proof since it is exactly the same with linear Schr\"odinger operator $e^{it\Delta}$.

Now, we introduce two Galilean operators, they are
\begin{equation}
J(t)=x\sinh t +i \cosh t \nabla_x, \quad  H(t)=x\cosh t+i\sinh t \nabla_x,
\label{8}
\end{equation}
conversely, $x$ and $\nabla_x$ can be expressed in terms of $J(t)$ and $H(t)$,
\begin{equation}
x=\cosh t H(t)-\sinh t J(t),\quad  i\nabla_x =\cosh t J(t)-\sinh t H(t).
\label{9}
\end{equation}
Furthermore, $J(t)$ and $H(t)$ enjoy the following property,

{\bf Lemma 2.3:} The operators $J$ and $H$ satisfy

\noindent 1. They are Heisenberg observables and consequently commute with the linear operator,
\begin{eqnarray*}
&& J(t)=U(t)i\nabla U(-t),\quad H(t)=U(t)xU(-t),\\
&& [i\partial_t+\frac {\Delta}2+\frac {|x|^2}2,J(t)]=[i\partial_t+\frac {\Delta}2+\frac{|x|^2}2, H(t)]=0.
\end{eqnarray*}
\noindent 2. \ They can be factorized as follows, for $t\neq 0$,
\begin{eqnarray*}
J(t)=i\cosh t e^{i\frac {|x|^2}2 \tanh t}\nabla_x(e^{-i\frac {|x|^2}2 \tanh t}\cdot),\\
H(t)=i\sinh t e^{i\frac {|x|^2}2 \coth t}\nabla_x(e^{-i\frac {|x|^2}2 \coth t}\cdot) .
\end{eqnarray*}
\noindent 3. Let $F\in C^1(\mathbf C,\mathbf C)$ and $F(z)=G(|z|^2)z$, then,
\begin{eqnarray*}
&&J(t)F(u)=\partial_z F(u) J(t)u-\partial_{\bar z}F(u)\overline{J(t)u},\\
&&H(t)F(u)=\partial_z F(u) H(t)u-\partial_{\bar z}F(u)\overline{H(t)u}.
\end{eqnarray*}
\noindent 4. There are embeddings(for instance),
\begin{eqnarray*}
\|f\|_{10}\le \|J(t)f\|_{\frac {30}{13}}, \quad \forall t\in \mathbf R,\\
\|f\|_{18}\le \|J(t)f\|_{\frac {18}{7}}, \quad \forall t\in \mathbf R.
\end{eqnarray*}
{\bf Proof:} The first point is easily checked thanks to (\ref{8}). The second one holds by direct computation, and
implies the last two one.

Formally, the solution of (\ref{1.1})-(\ref{1.2}) satisfies the following two conservation laws,
\begin{eqnarray*}
&& \mbox{Mass:} \quad M=\|u(t)\|_2=\|u_0\|_2,\\
&& \mbox{Energy:}\quad E(t)=\frac 12 \|\nabla u(t)\|_2^2-\frac 12
\|xu\|_2^2+\frac 13\|u(t)\|_6^6=const.
\end{eqnarray*}
As mentioned in the introduction, we split $E(t)$ by two ways.
First, define
$$
E_1(u(t))=\frac 12 \|\nabla u(t)\|_2^2+\frac 13
\|u(t)\|_6^6,\qquad E_2(u(t))=\frac 12 \|xu\|_2^2,
$$
it follows easily that,
$$
E(u(t))=E_1(u(t))+E_2(u(t)).
$$
Next, we define
\begin{eqnarray}
\begin{array}{ll}
\mathcal E_1(t):=\frac 12 \|J(t)u(t)\|_{L^2}^2+\frac 13\cosh^2 t \|u(t)\|_6^6,\\
\mathcal E_2(t):=\frac 12 \|H(t)u\|_{L^2}^2+\frac 13\sinh^2 t\|u(t)\|_6^6 ,
\end{array}
\label{8.1}
\end{eqnarray}
we see that $\mathcal E_1(t)$ and $\mathcal E_2(t)$ coincide with
$E_1(t)$ and $E_2(t)$ only at $t=0$. Furthermore, we have,

{\bf Lemma 2.4:} We can verify that:

1. $\mathcal E_1$ and $\mathcal E_2$ satisfy,
\begin{eqnarray}
&& \mathcal E_1(t)-\mathcal E_2(t)=E(t),\label{10}\\
&&\frac{d\mathcal E_1(t)}{dt}=\frac{d\mathcal E_2(t)}{dt}=-\frac23\sinh (2t)\|u(t)\|_6^6.\label{11}
\end{eqnarray}

2. The potential energy $\|u(t)\|_6^6 $ has exponentially decay in time:
\begin{eqnarray*}
\|u(t)\|_6^6\le 3E_1(0)\cosh^{-6}t, \quad\forall \ t\in\mathbf R.
\end{eqnarray*}

3. $\forall t\in \mathbf R$,
\begin{eqnarray}
&&\mathcal E_1(t)\le \mathcal E_1(0)=E_1(0),\label{12} \\
&&\|H(t)u(t)\|_2\le \|H(0) u(0)\|_2= \|xu_0\|_2=E_2(0).\label{13}
\end{eqnarray}
{\bf Proof:} The first point can be verified by (\ref{8}) and the equation (\ref{1.1}), see\cite{Carles1} for details.
Now let us
prove the second point. Integrating in time from $0$ to $t$, we see from (\ref{9}) that,
$$
\mathcal E_1(t)=\mathcal E_1(0)-\frac23\int_0^t\sinh (2s)\|u(s)\|_6^6ds.
$$
By (\ref{8.1}), we have
\begin{eqnarray*}
\cosh^2t\|u(t)\|_6^6&\le& 3\mathcal E_1(0)-2\int
\sinh(2s)\|u(s)\|_6^6ds\\
&=&3\mathcal E_1(0)-2\int_0^t\frac{\sinh(2s)}{\cosh^2s}\cosh^2s\|u(s)\|_6^6ds.
\end{eqnarray*}
Applying the Gronwall inequality yields:
$$
\cosh^2t\|u(t)\|_6^6\le 3\mathcal E_1(0)\exp\biggl[-2\int_0^t\frac{\sinh(2s)}{\cosh^2s}ds\biggr].
$$
Noting by direct computation,
$$
\int_0^t\frac{\sinh(2s)}{\cosh^2 s}ds=\ln \cosh t,
$$
thus, we have
$$
\cosh^2t\|u(t)\|_6^6\le 3\mathcal E_1(0)\cosh^{-4}t,
$$
and,
$$
\|u(t)\|_6^6\le 3\mathcal E_1(0)\cosh^{-6} t.
$$
Now, let's prove the third point. First, (\ref{12}) is easily
verified by using (\ref{11}). Next, noting (\ref{8.1}), (\ref{9})
and energy conservation, we see that,
\begin{eqnarray*}
\frac 12\|H(t)u(t)\|_2^2+\frac13\sinh^2t\|u(t)\|_6^6&=&\mathcal
E_1(t)-E(t)\\
&\le&\mathcal E_1(0)-E(0)\\
&=&\frac12\|\nabla u_0\|_2^2+\frac13\|u_0\|_6^6-(\frac12\|\nabla
u_0\|_2^2-\frac12\|xu_0\|_2^2+\frac13\|u_0\|_6^6)\\
&=&\frac12\|xu_0\|_2^2.
\end{eqnarray*}
Thus we get
$$
\|H(t)u(t)\|_2\le \|xu_0\|_2,
$$
which is exactly (\ref{13}).

Before ending this Section, we give the main theorems of this
paper.

{\bf {Theorem 1}:} Let $u_0\in \Sigma$ be radial, then the Cauchy
problem (\ref{1.1})-(\ref{1.2}) has a unique global solution in
$C(\mathbf R;\Sigma)\cap L^{10}_{xt}(\mathbf R\times \mathbf R^3)$
and satisfies
\begin{eqnarray}
\|u\|_{L^{10}(\mathbf R;L^{10})}\le C(\|\nabla
u_0\|_2,\|xu_0\|_2)\label{d},
\\
\max_{A\in\{J,H,I\}}\|A(\cdot) u\|_{L^q(\mathbf R;L^r)}\le
C(\|u_0\|_{\Sigma})\label{e}.
\end{eqnarray}
Furthermore, there exits a unique $u_+\in\Sigma$ such
that
$$
\|U(-t)u(t)-u_+\|_{\Sigma}\to 0,\mbox{ as } t\to \infty;
$$
there exists a unique $u_-\in\Sigma$ such that
$$
\|U(-t)u(t)-u_-\|_{\Sigma}\to 0, \mbox{ as } t\to -\infty.
$$

{\bf {Theorem 2}:}(Existence of wave operator)

Let $u_+\in\Sigma$ be radial, then there exists a unique solution
$u(x,t)$ of equation (\ref{1.1}) satisfying
\begin{eqnarray*}
&&\|u\|_{L^{10}(\mathbf R;L^{10})}\le C(E_1(u_+),E_2(u_+)),\\
&&\max_{A\in\{J,H,I\}}\|u\|_{L^{q}(\mathbf R;L^{r})}\le
C(\|u_+\|_{\Sigma}),\quad(q,r)\mbox{ admissible},
\end{eqnarray*}
and
$$
\|U(-t)u(t)-u_+\|_{\Sigma}\to0,\quad\mbox{as}\ t\to \infty;
$$
 Let $u_-\in\Sigma$ be radial, then there exists a unique
solution of equation (\ref{1.1}) satisfying
\begin{eqnarray*}
&&\|u\|_{L^{10}(\mathbf R;L^{10})}\le C(E_1(u_-),E_2(u_-)),\\
&&\max_{A\in\{J,H,I\}}\|u\|_{L^{q}(\mathbf R;L^{r})}\le
C(\|u_-\|_{\Sigma}),\quad(q,r)\mbox{ admissible},
\end{eqnarray*}
and
$$
\|U(-t)u(t)-u_-\|_{\Sigma}\to0,\quad\mbox{as}\ t\to -\infty.
$$

\section{Local wellposedness and global small solution}

In this section, we aim to get local solution and global small solution to the energy critical Schr\"odinger equation
with harmonic potential. By Duhamel's formula, it's enough to find solutions to the integral equation
$$
u(t)=U(t)u_0-i\int_0^tU(t-s)|u|^4u(s)ds.
$$
Define solution map by $\Phi(
u)(x,t)=U(t)u_0-i\int_0^tU(t-s)|u|^4u(s)ds$, one is required to
find fixed point of the map $\Phi$.

\bf{Proposition3.1} (Local wellposedness)

\rm For any $u_0\in\Sigma$, there exists maximal time interval $(T_*^-,T_*^+)$, such that (\ref{1.1})-(\ref{1.2}) has a
unique solution
$$
u(x,t)\in C((T_*^-,T_*^+);\Sigma)\cap L^{10}_{loc}((T_*^-,T_*^+);L^{10}).
$$
Furthermore, for any admissible pair $(q,r)$, one has
$$
\|A(\cdot)\|_{L^{q}_{loc}((T_*^-,T_*^+);L^{r})}<\infty, \quad A\in \{J,H,I\}.
$$
\bf{Proof:} \rm Let

$$
R=\max_{A\in\{J,H,I\}}\|A(\cdot)U(\cdot)u_0\|_
{L^{\frac{10}3}((T^-,T^+);L^{\frac{10}3})\cap L^{10}((T^-,T^+);L^{\frac{30}{13}})},
$$
with $T^-$, $T^+$ being specified later. By Strichartz estimate
Lemma2.2 and using the facts
$$
A(t)U(t)=U(t)B, \quad B\in \{i\nabla_x,x,I\},
$$
we see that $R$ is uniformly bounded w.r.t $T^-$, $T^+$,
and inparticulaly,
$$
R\le C\|u_0\|_{\Sigma},
$$
this in turn shows that $R$ is small when $T^+$ and $T^-$ are
small.

Define a set
$$
X=\biggl\{u(x,t)\biggl|
         \max_{A\in\{J,H,I\}}\|Au\|_{L^{\frac{10}3}((T^-,T^+);L^{\frac{10}3})\cap
         L^{10}((T^-,T^+);L^{\frac{30}{13}})}
             \le2R\biggr\},
$$
and the norm $\|\cdot\|_X$ is taken as the same as the one in the
capital bracket. First, we show that $X$ is stable under the
solution map $\Phi$. Choosing $u\in X$, and using Lemma2.2,
Lemma2.3, one computes that
$$
A(t)\Phi(u)(t)=A(t)U(t)u_0-i\int_0^t U(t-s)A(s)|u|^4u(s)ds,
$$
and,
\begin{eqnarray*}
&&\|A\Phi(u)\|_{L^{\frac{10}3}((T^-,T^+);L^{\frac{10}3})\cap L^{10}((T^-,T^+);L^{\frac{13}{30}})}\\
&=& R+\|\int_0^tU(t-s)A(s)|u|^4u(s)ds\|_{L^{\frac{10}3}((T^-,T^+);L^{\frac{10}3})\cap L^{10}((T^-,T^+);L^{\frac{30}{13}})}.
\end{eqnarray*}
By Strichartz estimate and Lemma2.3, the second term can be controlled by
\begin{eqnarray*}
&&C\|A(\cdot)|u|^4u\|_{L^{\frac{10}7}((T^-,T^+),L^{\frac{10}7})}\\
 &\le&
 C\|u\|_{10}^4\|A(\cdot)u\|_{L^{\frac{10}3}((T^-,T^+),L^{\frac{10}3})},
 \end{eqnarray*}
by embedding, this is smaller than $C\|u\|_X^5$. Consequently , we obtain
$$
\|A(\cdot)\Phi(u)\|_X\le R+C(2R)^5.
$$
Thus, $X$ is stable if $R$ is such that $C2^5R^4<\frac 14$. This
is available by choosing $T^-$ and $T^+$ small enough. Donate the
metric on $X$ by
$$
d(u_1,u_2)=\|u_1-u_2\|_{L^{\frac{10}3}((T^-,T^+);L^{\frac{10}3})},
$$
we need only to prove the contraction under this weak metric.
Taking $u_1$, $u_2$ $\in X$, we have
\begin{eqnarray*}
d(\Phi(u_1),\Phi(u_2))&=&\|\int_0^t
U(t-s)(|u_1|^4u_1(s)-|u_2|^4u_2)(s)ds\|_{L^{\frac{10}3}((T^-,T^+);L^{\frac{10}3})}\\
&\le &
C\||u_1|^4u_1-|u_2|^4u_2\|_{L^{\frac{10}7}((T^-,T^+);L^{\frac{10}7})}\\
&\le & C(\|u_1\|_{L^{10}((T^-,T^+);L^{10})}^4
+\|u_2\|_{L^{10}((T^-,T^+);L^{10})}^4)\|u_1-u_2\|_
       {L^{\frac{10}3}((T^-,T^+);L^{\frac{10}3})}\\
&\le & C(\|u_1\|_X^4+\|u_2\|_X^4)d(u_1,u_2)\\
&\le& 2C(2R)^4 d(u_1,u_2)\\
&\le&\frac 12d(u_1,u_2).
\end{eqnarray*}
Now, we now first fix $R$ such that $C2^5R^4<\frac14$, then choose
$T^-,T^+$ such that
$$
\max_{A\in\{J,H,I\}}\|A(\cdot)U(\cdot)u_0\|_ {L^{\frac{10}3}((T^-,T^+);L^{\frac{10}3})\cap
 L^{10}((T^-,T^+);L^{\frac{30}{13}})}<R,
 $$
by the fixed point theorem, we get a solution on $[T^-, T^+]$.
Once this is done, we extend this solution to the maximal time
interval $(T_*^-, T_*^+)$. The regularity property of the solution
follows from the Strichartz estimate. Thus, we conclude the proof
of Proposition3.1.

\bf{Proposition3.2}\rm (Global small solution)

 There exists an
absolute constant $\varepsilon>0$ such that when $u_0\in \Sigma$
and
$$
\|\nabla_x u_0\|_2\le \varepsilon,$$
(\ref{1.1})-(\ref{1.2}) has a unique global solution satisfying
\begin{eqnarray}
&&\|u\|_{L^{10}(\mathbf R;L^{10})}\le 2C\varepsilon,\label{a.1}\\
&&A(t)u(t,x)\in L^q(\mathbf R;L^r),\quad(q,r) \mbox{ admissible},
A\in\{J,H,I\}.\label{a}
\end{eqnarray}
Furthermore, there exists a unique function $u_+\in\Sigma$ such
that
$$
\|U(-t)u(t)-u_+\|_{\Sigma}\to 0,\quad\mbox{as } t\to \infty,
$$
and there exists a unique function $u_-\in\Sigma$ such that
$$
\|U(-t)u(t)-u_-\|_{\Sigma} \to 0,\quad\mbox{as } t\to-\infty.
$$
\bf{Proof:} \rm Being slightly different from the proof of Proprsition3.1, we define
\begin{eqnarray*}
R:&=&\|\nabla u_0\|_2,\\
X:&=&\biggl\{u(x,t)\biggl| \quad \|u\|_{L^{\frac{10}3}(\mathbf R;L^{\frac{10}3})}\le
2C\|u_0\|_2,\\
&&\qquad\qquad\qquad \|H(\cdot)u\|_{L^{\frac{10}3}(\mathbf R;L^{\frac{10}3})}\le
2C\|xu_0\|_2,\\
&&\qquad\qquad\qquad \|J(\cdot)u\|_{L^{\frac{10}3}(\mathbf R;L^{\frac{10}3})\cap L^{10}(\mathbf R;L^{\frac{30}{13}})}
\le 2CR\biggr\},
\end{eqnarray*}
where, $C$ in the bracket is the Strichartz constant. First, we show $\Phi$ is onto from $X$ to $X$. Taking $u\in X$,
we verify that
\begin{eqnarray*}
&&\|J(\cdot)\Phi(u)\|_{L^{\frac{10}3}(\mathbf R;L^{\frac{10}3})\cap L^{10}(\mathbf
R;L^{\frac{30}{13}})}\\
&\le& C\|\nabla u_0\|_2+C\|J(\cdot)|u|^4u\|_{L^{\frac{10}7}(\mathbf
R;L^{\frac{10}7})}\\
&\le& C R+C\|u\|_{L^{10}(\mathbf R;L^{10})}^4\|J(\cdot)u\|_ {L^{\frac{10}3}(\mathbf
R;L^{\frac{10}3})}\\
&\le& CR+C\|J(\cdot)u\|^5_{L^{\frac{10}3}(\mathbf R;L^{\frac{10}3})\cap L^{10}(\mathbf R;L^{\frac{30}{13}})}\\
&\le& CR+C(2CR)^5.
\end{eqnarray*}
Here, we have used Strichartz estimate and embedding. If $R$ is taken such that
\begin{equation}
C(2CR)^5<\frac 12CR,
\label{b}
\end{equation}
then we obtain
$$
\|J(\cdot)\Phi(u)\|_{L^{\frac{10}3}(\mathbf R;L^{\frac{10}3})\cap L^{10}(\mathbf R;L^{\frac{30}{13}})}\le 2CR.
$$
Now, we verify the first two properties in the bracket. Taking $u\in X$, we have
\begin{eqnarray*}
\|\Phi(u)\|_{L^{\frac{10}3}(\mathbf R;L^{\frac{10}3})}&\le &C\|u_0\|_2+C\||u|^4u\|_
{L^{\frac{10}7}(\mathbf R;L^{\frac{10}7})}\\
&\le& C\|u_0\|_2+C\|J(\cdot)u\|_ {L^{10}(\mathbf R;L^{\frac{30}{13}})}^4\|u\|_{L^{\frac{10}3}(\mathbf
R;L^{\frac{10}3})}\\
&\le& C\|u_0\|_2+C(2CR)^4 2C\|u_0\|_2\\
&\le& 2C\|u_0\|_2.
\end{eqnarray*}
By the same way, one gets
$$
\|H(\cdot)\Phi u\|_{L^{\frac{10}3}(\mathbf R;L^{\frac{10}3})}\le 2C\|xu_0\|_2.
$$
Hence, $\Phi$ is a map from $X$ to $X$. To complete the proof of
the existence part, we donate $X$ with the weak metric
$$
d(u_1,u_2)=\|u_1-u_2\|_{L^{\frac{10}3}(\mathbf R;L^{\frac{10}3})},
$$
and plan to prove the contraction under this metric. Choosing
$u_1$, $u_2$ $\in X$ with same data, we have
$$
\Phi(u_1)(t)-\Phi(u_2)(t)=-i\int_0^t U(t-s)(|u_1|^4u_1-|u_2|^4u_2)(s)ds.
$$
Taking Lebesgues norm on each sides of the equation and using Strichartz estimate and embedding, we easily get
$$
d(\Phi(u_1),\Phi(u_2))\le \frac12 d(u_1,u_2).
$$
By the fixed point theorem, we obtain a unique solution $u\in X$. The properties (\ref{a}) and (\ref{a.1})
follows directly from the Strichartz estimate and embedding.
Now, Let's turn to the proof of the second part. Let $u_+(x,t)=i\int_0^{\infty}U(-s)|u|^4u(x,s)ds$, one sees that
$$
U(-t)u(t)-u_+(t)
=i\int_t^{\infty} U(-s)|u|^4u(s)ds.
$$
Noting that by Lemma2.3, there holds
\begin{eqnarray*}
i\nabla_x(U(-t)u(t)-u_+(t))&=&i\int_t^{\infty}
U(-s)J(s)|u|^4u(s)ds,\\
x(U(-t)u(t)-u_+(t))&=&i\int_t^{\infty}
U(-s)H(s)|u|^4u(s)ds,
\end{eqnarray*}
We see that
$$
\|U(-t)u(t)-u_+(t)\|_{\Sigma}\le \max_{A\in\{J,H,I\}}\|\int_t^{\infty} U(-s)A(s)|u|^4u(s)ds\|_2.
$$
Noting that the operator $U(\cdot)$ is unitary on $L^2$, we can bound the right side of the above equation by
$$
\|u\|_{L^{10}((t,\infty);L^{10})}^4\|A(\cdot)u\|_{L^{\frac{10}3}((t,\infty);L^{\frac{10}3})}
$$
which tends to 0 as $t$ tends to $\infty$. The scattering in the negative direction follows from the way.
This finally gives Proposition3.2.

The above two proposition provide no answer about whether the solution with large data is global. Assume $u$
is the local solution on the maximal time interval $(T_*^-,T_*^+)$, our first purpose is to show that $T_*^-=-\infty$,
and $T_*^+=\infty$ once we prove an apriori estimate on $u$ as follows:

\bf{Lemma3.3} \rm Assume $u$ be a maximal solution on $(T_*^-,T_*^+)$ with finite energy. Then if for any $I
\in (T_*^-,T_*^+)$, $u$ satisfies
\begin{equation}
\|u\|_{L^{10}(I;L^{10})}<C(|I|, E_1(u_0),E_2(u_0)),
\label{j.10}
\end{equation}
then $T_*^-=-\infty$, $T_*^+=\infty$.
Here, $|I|$ denotes the length of I.

\bf{Proof:} \rm Let's discuss in the positive time direction. Assume otherwise that $T_*^+<\infty$, we will get
a contradiction by showing the solution can be extended beyond $T_*^+$. Our strategy is as follows: we first take
$t_0\in[0, T_*^+)$ that is close enough to $T_*^+$, then aim to solve the same Cauchy problem from $t_0$ forward. Once
we have shown that there exists $\delta>0$ such that
\begin{equation}
\max_{A\in\{J,H,I\}}\|A(\cdot)U(\cdot-t_0)u(t_0)\|_
{L^{\frac{10}3}([t_0,T_*^++\delta);L^{\frac{10}3})\cap L^{10}([t_0,T_*^++\delta);L^{\frac{30}{13}})}
<R,
\label{j.11}
\end{equation}
where $R$ is a same constant in Proposition3.1, we will establish a contraction mapping. This allows us to
extend the solution at least beyond $T_*^++\delta$ and contradicts the maximum property of $T_*^+$. So, let's prove
(\ref{j.11}). First of all, by Strichartz estimate and some routine arguments,
we see that (\ref{j.10}) implies that
\begin{equation}
\max_{A\in \{J,H,I\}}\|A(\cdot)u\|_{L^q([0, T_*^+);L^r)}\le C(|T^+_*|,\|u_0\|_{\Sigma}),\ (q,r) \mbox{ admissible}.
\label{j.12}
\end{equation}
By Duhamel's formula, we see that $u$ solves the equation
$$
u(t)=U(t-t_0)u(t_0)-i\int_{t_0}^t U(t-s)|u|^4u(s)ds,\mbox{ for } t\in [t_0, T_*^+),
$$
and hence,
\begin{eqnarray*}
A(t)U(t-t_0)u(t_0)&=&A(t)u(t)+i\int_{t_0}^t
A(t)U(t-s)|u|^4u(s)ds\\
&=& A(t)u(t)+i\int_{t_0}^t U(t-s)A(s)|u|^4u(s)ds.
\end{eqnarray*}
Taking Lebesgues norm on each side to the equation, one gets
\begin{eqnarray}
&&\|A(\cdot)U(\cdot-t_0)\|_{L^{\frac{10}3}([t_0,T_*^+);L^{\frac{10}3})\cap
L^{10}([t_0,T_*^+);L^{\frac{30}{13}})}\nonumber\\
&&\le \|A(\cdot)u\|_{L^{\frac{10}3}([t_0,T_*^+);L^{\frac{10}3})\cap L^{10}([t_0,T_*^+);L^{\frac{30}{13}})}
 +C\|A(\cdot)|u|^4u\|_{L^{\frac{10}7}([t_0,T_*^+);L^{\frac{10}7})}.
 \nonumber \\
 \label{j.13}
\end{eqnarray}
Having (\ref{j.12}) in mind and applying H\"older, we see (\ref{j.13}) is smaller than $\frac R2$ if $t_0$ is close enough to
 $T_*^+$.

On the other hand, since by the Strichartz estimates
\begin{eqnarray*}
&&\qquad\|A(\cdot)U(\cdot-t_0)\|_{L^{\frac{10}3}(R;L^{\frac{10}3})\cap
L^{10}(\mathbf R;L^{\frac{30}{13}})}\\
&&=\|U(\cdot - t_0)
A(t_0)u(t_0)\|_{L^{\frac{10}3}(R;L^{\frac{10}3})\cap
L^{10}(\mathbf
R;L^{\frac{30}{13}})}\\
&&\le C\|A(t_0)u(t_0)\|_2 \le C(\|u_0\|_{\Sigma}),
\end{eqnarray*}
thus, once $t_0$ is fixed by $(\ref{j.13})\le \frac R2$, one is
allowed to choose a $\delta>0$ sufficiently small such that
\begin{equation}
\|A(\cdot)U(\cdot-t_0)\|_{L^{\frac{10}3}([T_*^+,T_*^++\delta);L^{\frac{10}3})\cap
L^{10}([T_*^+,T_*^++\delta);L^{\frac{30}{13}})} \le \frac R2 ,
\label{14}
\end{equation}
(\ref{j.11}) then follows by collecting $(\ref{j.13})\le \frac R2$
and (\ref{14}), also implies Lemma3.3.

Lemma3.3 says the solution is global in the sense that it exists on arbitrary finite time interval $(-T,T)$. Inparticularly
, it doesn't imply that the solution enjoy certain global space-time estimate which is the usual requirement in the scattering
theory. However, we can complement this by the decay estimate Lemma2.4.

{\bf Lemma3.4:} Assume $u$ be the global solution in the sense above, then $u$ satisfies
\begin{equation}
\max_{A\in \{J,H,I\}} \|A(\cdot)u\|_{L^q(\mathbf R;L^r)}\le C(\|u_0\|_{\Sigma}),\quad \forall (q,r)\quad \mbox{admissible}.
\label{14.1}
\end{equation}
and there is scattering.

{\bf Proof:} Fixing a small number $\varepsilon$ and taking $T\ge
T_0=(\frac {3E_1(0)}{\varepsilon^6})^{\frac 16}$, we have
$$
3E_1(0) \cosh^{-6} T \le \varepsilon^6 ,
$$
thus by the decay estimate Lemma2.4, one has
\begin{equation}
\|u\|_{L^{\infty}([T,\infty);L^6)}\le \varepsilon .
\label{14.2}
\end{equation}
By Duhamel's formula,  on $[T,\infty)$, $u$ satisfies the equation
$$
u(t)=U(t-T)u(T)-i\int_T^tU(t-s)|u|^4 u(s) ds,
$$
taking a special admissible pair $(6,\frac {18}7)$ and applying Strichartz estimate gives,
\begin{eqnarray*}
&&\|J(\cdot)u\|_{L^6([T,\infty);L^{\frac {18}7})}\\
&\le& \|J(\cdot) U(\cdot-T)u(T)\|_{L^6([T,\infty);L^{\frac {18}7})}+
\|\int_T^t J(\cdot)U(\cdot-s)|u|^4u(s)ds\|_{L^6([T,\infty);L^{\frac
{18}7})}\\
&\le& C\|J(T)u(T)\|_2+C\|J(\cdot)|u|^4u\|_{L^{\frac 32}([T,\infty);L^{\frac {18}{13}})}
\end{eqnarray*}
The first term is smaller than $CE_1(u_0)^{\frac12}$, and by H\"older, the second term is controlled by
$$
C\|u\|_{L^{\infty}([T,\infty);L^6)}\|u\|^3_{L^6([T,\infty);L^{18})}\|J(\cdot)u\|_{L^6([T,\infty);L^{\frac
{18}7})},
$$
in view of (\ref{14.2}) and embedding, we further estimate it by
$$
C\varepsilon \|J(\cdot)u\|^4_{L^6([T,\infty);L^{\frac {18}7})}.
$$
Hence we get an estimate for $J(\cdot)u$ as follows,
$$
\|J(\cdot)u\|_{L^6([T,\infty);L^{\frac {18}7})}\le C E_1(u_0)^{\frac 12}+C\varepsilon\|J(\cdot)u\|^4_{L^6([T,\infty);
L^{\frac {18}7})}.
$$
(The more rigorous way is to do estimate on $[T,R], \quad R<\infty$, then take supreme w.r.t. $R$.) This implies that
$\|J(\cdot)u\|_{L^6([T,\infty);L^{\frac {18}7})}$ is bounded if $\varepsilon$ is smaller than a constant which depends
only on $E_1(u_0)$. Once this has been obtained, one can get
\begin{eqnarray*}
\|J(\cdot)u\|_{L^6([0,\infty);L^{\frac {18}7})} &\le& \|J(\cdot)u\|_{L^6([0,T_0);L^{\frac {18}7})}+
\|J(\cdot)u\|_{L^6([T_0,\infty);L^{\frac {18}7})}\\
&\le& C(E_1(u_0),E_2(u_0)).
\end{eqnarray*}
By time reversing and Strichartz estimate, we obtain (\ref{14.1}).

Having Lemma3.3 and Lemma3.4 in mind, in order to prove Theorem1.1, we need only to show
$$
\|u\|_{L^{10}([-T_0,T_0];L^{10})}\le C(E_1(u_0),E_2(u_0)),
$$
where $T_0$ is defined in Lemma3.4 and depends only on $E_1(u_0)$.

Now, we fix two constants $E$ and $B$ such that
$$
E_1(u_0)=E, \qquad  E_2(u_0)=B ,
$$
then our task becomes to prove
\begin{equation}
\|u\|_{L^{10}([-T_0(E),T_0(E)];L^{10})}\le C(E,B),
\label{c}
\end{equation}
if $u$ is a finite energy solution on $[-T_0(E),T_0(E)]$ with
$u(0)=u_0$. From (\ref{9}) and Lemma2.4, we compute that,
\begin{eqnarray*}
\|\nabla u(t)\|_2 &\le& \cosh t \|H(t)u(t)\|_2+\sinh t
\|J(t)u(t)\|_2\\
&\le& C(E,B),\qquad \forall t \in [-T_0(E),T_0(E)],\\
\|xu(t)\|_2&\le& \sinh t\|H(t)u(t)\|_2+\cosh t\|J(t)u(t)\|_2\\
&\le& C(E,B),\qquad \forall t \in [-T_0(E),T_0(E)],
\end{eqnarray*}
Thus, there exists $\Lambda_1(E,B)$ and $\Lambda_2(E,B)$ such that
$$
E_1(u(t))\le \Lambda_1, \quad E_2(u(t))\le \Lambda_2 ,\qquad
\forall t\in\ [-T_0(E),T_0(E)].
$$
If there exists $\eta_1=\eta_1(\Lambda_1,\Lambda_2)$ such that on
every time interval $I$ with length $2\eta_1^4$, one has
$$
\|u\|_{L^{10}(I;L^{10})}\le C(\Lambda_1,\Lambda_2),
$$
then we can divide $[-T(E),T(E)]$ into $O(\frac {T(E)}{\eta_1^4})
$ subintervals, and get (\ref{c}) by summing the estimates on each
subinterval. Indeed, we plan to prove the following proposition.

{\bf Proposition3.5:} Let $t'\in \mathbf R$ be arbitrarily fixed and $u(t')\in \Sigma$ satisfying
$$
E_1(u(t'))\le \Lambda_1,\qquad E_2(u(t'))\le \Lambda_2,
$$
then we have a small constant $\eta_1$ which depends only on $
(\Lambda_1,\Lambda_2+\Lambda_1) $ such that
the Cauchy problem of (\ref{1.1}) with prescribed data $u(t')$ at time $t'$ is at least solvable on
$[t'-\eta_1^4,t'+\eta_1^4] $ and the solution $u$ satisfy
$$
\|u\|_{L^{10}([t'-\eta_1^4,t'+\eta_1^4];L^{10})}\le
C(\Lambda_1,\Lambda_2).
$$

Assuming Proposition3.5 hold true, let's give a remark about the
proof of Theorem1, Theorem2. First of all, in Theorem1, we are
left to prove the regularity part and the scattering
part(\ref{d}), (\ref{e}) of the global solution which can be
deduced from the Strichartz estimate and some routine arguments.
see\cite{Bourgain} and the proof of Proposition3.2 for details.
The proof Theorem2 is a bit different, so we sketch it below.

\bf{Proof of Theorem2:} \rm We need only to show the integral equation
\begin{equation}
u(t)=U(t)u_++i\int_t^{\infty} U(t-s)|u|^4u(s)ds,\label{f}
\end{equation}
has a unique global solution with global spacetime estimates.
First of all, we seek for local solution. Define the solution map
by $\Phi(u)(t)= U(t)u_++i\int_t^{\infty} U(t-s)|u|^4u(s)ds$, and
denote $R=\|\nabla u_+\|_2$. By choosing $T=T(R)$ large enough,
say, $\cosh T\ge CR$, we see that $\Phi$ is a contraction map on
the set
\begin{eqnarray*}
X=\biggl\{u(x,t);\quad
\begin{array}{l}
\|u\|_{L^{\frac{10}3}([T,\infty);L^{\frac{10}3})}\le
2C\|u_+\|_2,\\
\|H(\cdot) u\|_{L^{\frac{10}3}([T,\infty);L^{\frac{10}3})}\le
2C\|xu_+\|_2,\\
\|J(\cdot )u\|_{L^{\frac{10}3}([T,\infty);L^{\frac{10}3})\cap
L^{10}([T,\infty);L^{\frac{30}{13}})}\le 2CR
\end{array}
\biggr\},
\end{eqnarray*}
donated with the metric
$d(u_1,u_2)=\|u_1-u_2\|_{L^{\frac{10}3}([T,\infty);L^{\frac{10}3})}$.
The proof is routine, except needing to notifying that the gain
$\cosh^{-1}T$ from the embedding
$\|u\|_{L^{10}([T,\infty);L^{10})}\le
C\cosh^{-1}T\|J(\cdot)u\|_{L^{10}([T,\infty);L^{\frac{30}{13}})}$
gives the dependence of $T$ on $R$. Once we get the local
solution, we can find a finite time $T=T(E_1(u_+))$ such that
$u(T)\in\Sigma$ and has the bound only depending on
$\|u_+\|_{\Sigma}$. At that moment, by solving a finite time
Cauchy problem, we are allowed to get global solution of (\ref{f})
by Theorem1 and the uniqueness. The scattering part of Theorem2
follows easily from the global spacetime bound of $u(x,t)$, thus
we end the proof of Theorem2.

The remaining part of the paper is devoted to the proof of
Proposition3.5, by time translation, we may assume $t'=0$. We
begin the proof by giving the Morawetz estimate of the
Schr\"odinger equation with potential in the following section.

\section{Morawetz estimate for solutions of Schr\"odinger
equations with potential}

We first give the local mass conservation of $u$. Taking a smooth
function $\chi(x)\in C_0^{\infty}(\mathbf R^3)$ such that $\chi(x)=1$ if $|x|\le \frac 12$ and $\chi(x)=0$ if
$|x|\ge 1$. Then we have,

{\bf Proposition4.1:} Let $u$ be the smooth solution of (\ref{1.1}), and define local mass of $u$ to be
$$
Mass(u(t),B(x_0,R))=(\int \chi^2(\frac {x-x_0}R)|u(t,x)|^2dx)^{\frac 12} ,
$$
then,
\begin{equation}
\partial_t Mass(u(t),B(x_0,R))\le \frac {\|u(t)\|_{\dot H^1}}R,
\label{17}
\end{equation}
and
\begin{equation}
Mass(u(t),B(x_0,R))\le R\|u(t)\|_{\dot H^1} .
\label{18}
\end{equation}
For the self-containedness of the paper, we give the proof of this proposition.

{\bf Proof:} Noting that $u$ satisfies the equation (\ref{1.1}),
we have
$$
\partial_t Mass(u(t),B(x_0,R))^2=\int \chi^2(\frac {x-x_0}R) 2Re\biggl[\bar
u (\frac i2 \Delta u +\frac {i}2|x|^2 u-i|u|^4 u)\biggr](x)dx,
$$
which equals to
$$
-\int \chi^2(\frac {x-x_0}R)Im(\bar u\Delta u)(x)dx ,
$$
by a simple computation. Using integrating by parts, we finally get
$$
\partial_t Mass(u(t),B(x_0,R))^2=\frac 2R \int \chi(\frac
{x-x_0}R) (\nabla \chi)(\frac {x-x_0}R)Im(u\nabla \bar u)(x,t)dx.
$$
By H\"older inequality, the right hand side can be controlled by
$$
\frac 2R Mass(u(t),B(x_0,R))\|\nabla u(t)\|_2,
$$
from which (\ref{17}) follows. Now, let's prove (\ref{18}). Using Hardy's inequality, one has
\begin{eqnarray*}
Mass(u(t),B(x_0,R))^2&=& \int \chi^2(\frac {x-x_0}R)|u(t,x)|^2
dx\\
&\le&\sup_{x\in \mathbf R^3} \chi^2(\frac {x-x_0}R) |x-x_0|^2 \int \frac {|u(t,x)|^2}{|x-x_0|^2}
dx\\
&\le & R^2\|u(t)\|^2_{\dot H^1}.
\end{eqnarray*}
Thus, we get (\ref{18}).

{\bf Proposition4.2:}\ (Morawetz inequality) \ Let $u$ be the solution of (\ref{1.1}) with finite energy. Then we have
\begin{eqnarray}
&&\int_I \int_{|x|\le A|I|^{\frac 12}} \frac {|u(t,x)|^6}{|x|} dx
dt \nonumber \\
&\le& CA|I|^{\frac 12}
(\|\nabla u\|^2_{L^{\infty}(I;L^2)}+\|xu\|^2_{L^{\infty}(I;L^2)}+\|u\|^6_{L^{\infty}(I;L^6)})
\mbox {  for all  } A\ge 1  .                \nonumber \\
\label{19}
\end{eqnarray}

{\bf Proof:} We prove this result by following the idea in \cite{Tao}.
 Assume without loss of generality that $u$ is a smooth solution of (\ref{1.1}). First, by a direct computation,
we get
\begin{equation}
\partial_t Im(\partial_k u\bar u)=\frac 14 \partial_k \Delta
(|u|^2)-Re \partial_j(\bar u_k u_j)-\frac 23 \partial_k(|u|^6)+x_k|u|^2,
\label{20}
\end{equation}
here, we use $\partial_k f$ or $f_k$ to denote $\frac {\partial f(x)}{\partial x_k}$. Let $a(x)$ be a smooth radial solution
to be choose later. Multiplying (\ref{20}) by $a_k(x)$ and integrating on $\mathbf R^3$, we get
\begin{eqnarray}
&&\partial_t\int_{\mathbf R^3} Im(\partial_k u\bar u)(x) a_k(x)dx=
\int a_{jk}(x) Re(\bar u_k u_j)(x)dx \nonumber\\
&&-\frac 14 \int \Delta \Delta
a(x)|u|^2(x)dx
+\frac 23\int \Delta a(x)|u|^6(x)dx+\int a_k(x)x_k|u|^2(x)dx.\nonumber \\\label{21}
\end{eqnarray}
Taking $\chi(x)\in C_0^{\infty} (\mathbf R^3)$ satisfying $\chi(x)=1$ as $|x|\le 1$ and $\chi(x)=0$ as $|x|\le 2$. Letting
$a(x)=(\varepsilon^2+|x|^2)^{\frac 12}\chi(\frac xR)$, we claim that, on $|x|\le R$,
\begin{eqnarray*}
&&a(x)=(\varepsilon^2+|x|^2)^{\frac 12}, \quad a_k(x)=\frac {x_k}{(\varepsilon^2+|x|^2)^{\frac
12}},\\
&&\Delta a(x)=\frac 2{(\varepsilon^2+|x|^2)^{\frac 12}}+\frac {\varepsilon^2}{(\varepsilon^2+|x|^2)^{\frac
32}},\quad \Delta \Delta a(x)=-\frac {15\varepsilon^2}{(\varepsilon^2+|x|^2)^{\frac
72}},\\
&&a_{jk}(x)Re(\bar u_k u_j)(x)\ge 0, \quad a_k(x)x_k=\frac {|x|^2}{(\varepsilon^2+|x|^2)^{\frac 12}}\ge 0.
\end{eqnarray*}
The first four points follow by directly differentiating $a(x)$ on $|x|\le R$. To see the fifth point, we do further computation,
\begin{eqnarray}
a_{jk}(x)Re(\bar u_k u_j)(x) &=& (\frac {\delta_{jk}}{(\varepsilon^2+|x|^2)^{\frac 12}} -
\frac {x_jx_k}{(\varepsilon^2+|x|^2)^{\frac 32}}) Re(\bar u_k u_j)(x)
\nonumber\\
&=& \frac {|\nabla u|^2}{(\varepsilon^2+|x|^2)^{\frac 12}}-\frac {Re(x_j u_j x_k\bar u_k)}{(\varepsilon^2+|x|^2)^{\frac
32}}\nonumber\\
&=&  \frac {|\nabla u|^2}{(\varepsilon^2+|x|^2)^{\frac 12}}-\frac{|x|^2|u_r|^2}
{(\varepsilon^2+|x|^2)^{\frac 32}} ,\label{22}
\end{eqnarray}
where $u_r$ denotes the radial derivative. Since $|u_r|\le|\nabla u|$, one sees that (\ref{22}) is greater than
$$
\frac {|\nabla u|^2\varepsilon^2}{(\varepsilon^2+|x|^2)^\frac 32}\ge 0 ,
$$
from which the positivity follows. The last point is an easy
consequence of the second one. Keeping the above claims in mind, we get
from (\ref{21}) that
\begin{eqnarray*}
&&\frac 43 \int_{|x|\le R}\frac{|u|^6}{(\varepsilon^2+|x|^2)^{\frac
12}}dx\\
&&\le \partial_t\int_{\mathbf R^3}Im(\partial_k u\bar u)(x) a_k(x)dx +\int_{R\le |x|\le 2R}
|a_{jk}(x)Re(\bar u_k u_j)(x)|\\&&+\frac 14|\Delta \Delta a(x)||u|^2(x)
+\frac 23 |\Delta a(x)||u|^6(x)+|a_k(x) x_k| |u|^2(x)dx .
\end{eqnarray*}
Integrating in time on $I$, we get
\begin{eqnarray}
&&\frac 43 \int_I\int_{|x|\le R} \frac
{|u|^6}{(\varepsilon^2+|x|^2)^{\frac 12}} dx \le
\sup_{t\in I} \biggl|\int _{\mathbf R^3} Im(\bar u u_k)(x)a_k(x)dx\biggr| \nonumber \\
&&+|I|\sup_{t\in I} \biggl(\int_{R\le |x|\le 2R} |a_{jk}(x) u_k
u_j(x,t)|+\frac 14|\Delta \Delta a(x)||u|^2(x,t)|\nonumber\\
&&\qquad\quad+
\frac 23|\Delta a(x)| |u|^6(x,t)|+|a_k(x) x_k| |u|^2(x,t)|dx\biggr).\nonumber \\
\label{23}
\end{eqnarray}
To estimate each term on the right hand side of the inequality, we
need rough bounds on the derivatives of $a$ as $R\le |x|\le 2R$,
they are
$$
\begin{array}{cc}
 \mbox{when}\  R\le |x|\le 2R; & |a_k(x)|\le C\frac {(\varepsilon^2+R^2)^{\frac 12}}R, \\
  & |a_{jk}(x)|\le C\frac {(\varepsilon^2+R^2)^{\frac 12}}{R^2}, \\
& |\Delta \Delta a(x)|\le C\frac {(\varepsilon^2+R^2)^{\frac
12}}{R^4}.
\end{array}
$$
Using these bounds, we can further control (\ref{23}) by
\begin{eqnarray*}
&&C(\varepsilon^2+R^2)^{\frac 12} \|\nabla
u\|^2_{L^{\infty}(I;L^2)}
+C|I|\frac {(\varepsilon^2+R^2)^{\frac
12}}{R^2}
\|\nabla u\|^2_{L^{\infty}(I;L^2)}\\ &&\qquad+C|I|\frac {(\varepsilon^2+R^2)^{\frac 12}}{R^2} \|u\|_{L^{\infty}(I;L^6)}^6
 +C|I|\frac {(\varepsilon^2+R^2)^{\frac 12}}{R^2} \|x u\|^2_{L^{\infty}(I;L^2)},
\end{eqnarray*}
which is smaller than
$$
C\biggl((\varepsilon^2+R^2)^{\frac 12}+|I|\frac
{(\varepsilon^2+R^2)^{\frac 12}}{R^2}\biggr) \biggl(\|\nabla u\|
^2_{L^{\infty}(I;L^2)}+\|xu\|^2_{L^{\infty}(I;L^2)}+\|u\|^6_{L^{\infty}(I;L^6)}\biggr).
$$
Choosing $R=A|I|^{\frac 12}$ and letting $\varepsilon \rightarrow
0$, (\ref{23}) becomes
\begin{eqnarray*}
&& \int_{I} \int_{|x|\le A|I|^{\frac 12}} \frac {|u|^6}{|x|}dx dt\le
C(A|I|^{\frac 12}+A^{-1}|I|^{\frac 12})(\|\nabla
u\|_{L^{\infty}(I;L^2)}^2+\|xu\|_{L^{\infty}(I;L^2)}^2+
 \|u\|_{L^{\infty}(I;L^2)}^6) \\
&&\le CA|I|^{\frac 12}(\|\nabla
u\|_{L^{\infty}(I;L^2)}^2+\|xu\|_{L^{\infty}(I;L^2)}^2+
\|u\|_{L^{\infty}(I;L^2)}^6),
\end{eqnarray*}
since $A\ge 1$. This is exactly (\ref{19}).

Being different from the Morawetz estimate for equations without
potential, the term
$$
\|\nabla u\|_{L^{\infty}(I;L^2)}^2+
 \|xu\|_{L^{\infty}(I;L^2)}^2+\|u\|_{L^{\infty}(I;L^6)}^6
 $$
cann't be substituted by a quantity independent on $I$ since it is
not conserved in time. However, because we have restricted this
problem in a finite time interval, we are allowed to control this
term. Indeed, we have

{\bf Corollary4.3:} Let $u$ be a finite energy solution on $I$ and satisfy
$$
E_1(u(t))\le C_1;\quad E_2(u(t))\le C_2,  \forall \ t\in I,
$$
then we have
\begin{equation}
\int_I\int_{|x|\le A|I|^{\frac 12}} \frac {|u(t,x)| ^6}{|x|}dxdt
\le C(C_1,C_2)A|I|^{\frac12}\ \mbox {  for all  } A\ge 1.
\label{24}
\end{equation}

\section{Paraproduct decomposition and Littlewood-Paley}

Keeping Corollary4.3 in mind, we begin to prove Proposition3.5
from this Section to the end. Thanks to the local solution theory,
we may assume the solution has been existed on
$[-\eta_1^4,\eta_1^4]$ and only aim to show the spacetime bound on
it. Let $\eta_1$ be a small number that meets all the conditions
that will appear in the proof, then dividing $[0,\eta_1^4]$ into
$J_1$ subintervals and $[-\eta_1^4,0]$ into $J_2$ subintervals
such that on each subinterval $I_j$, we have
$\eta_1\le\|u\|_{L^{10}(I_j,L^{10})}\le 2\eta_1$. So we are left
to control $J_1$, $J_2$ by constant $C(\Lambda_1,\Lambda_2)$.
Without loss of generality, we only do analysis in the positive
time direction. Following Bourgain\cite{Bourgain}, we classify the
subintervals into three components $I^{(1)}$, $I^{(2)}$,
$I^{(3)}$, and each contains $\frac{J_1}3$ consecutive
subintervals. It's on the middle component that we do most
analysis. Our first aim is to show the existence of a sequence of
bubbles somewhere in space at a sequence of times $t_j$ which
belongs to the subinterval $I_j$. We realize this by doing
analysis on one specified subinterval. At First, we show some
regularity property. (From this section to the end, the constant
$C$ may depend on $\Lambda_1$, $\Lambda_2$.)

{\bf Proposition 5.1:} Let $I_j$ be one of the subintervals, that is
 $I_j\subset [0,\eta_1^4]$ and
$$
\eta_1\le \|u\|_{L^{10}(I_j;L^{10})}\le 2\eta_1.
$$
Then $u$ satisfies
$$
\|\nabla u\|_{L^{\frac {10}3}(I_j;L^{\frac
{10}3})}+\|xu\|_{L^{\frac {10}3}(I_j;L^{\frac {10}3})}\le
C(\Lambda_1,\Lambda_2).
$$
{\bf Proof:} Noting that $I_j\subset [0,\eta_1^4]$, it's suffices
to prove the same space time bound for $J(t)u$ and $H(t)u$. By
Duhamel, on $I_j=[t_j,t_{j+1}]$, $u$ satisfies
$$
u(t)=U(t-t_j)u(t_j)-i\int_{t_j}^t U(t-s)|u|^4 u(s) ds.
$$
Let $A(t)\in \{ J(t),H(t)\}$, then we have
$$
A(t) u(t)=U(t-t_j)A(t_j)u(t_j)-i\int_{t_j}^t U(t-s)A(s)|u|^4
u(s)ds.
$$
Using Strichartz estimate, we get
\begin{eqnarray*}
\|A(\cdot) u\|_{L^{\frac {10}3}(I_j;L^{\frac {10}3})}&\le&
C\|A(t_j)u(t_j)\|_2 +C\|u\|^4_{L^{10}(I_j;L^{10})}\|A(\cdot)
u\|_{L^{\frac {10}3}(I_j;L^{\frac{10}3})}\\
&\le& C\|A(t_j)u(t_j)\|_2+C\eta_1^4\|A(\cdot)u\|_{L^{\frac {10}3}(I_j;L^{\frac{10}3})}.
\end{eqnarray*}
Noting that $\|A(t_j)u(t_j)\|_2$ is bounded by $C(\Lambda_1,\Lambda_2)$, we get
$$
\|A(\cdot)u\|_{L^{\frac {10}3}(I_j;L^{\frac {10}3})}\le C(\Lambda_1,\Lambda_2),
$$
by $C\eta_1^4<\frac 12$. Thus we end this Proposition.

{\bf Proposition5.2} Let $I_j$ be one of the subintervals, that is
$I_j\in[0,\eta_1^4]$ and $\eta_1\le \|u\|_{L^{10}(I_j;L^{10})}\le
2\eta_1$. Then there exists $t_j\in I_j$, $x_j\in \mathbf R^3$ and
$N\ge N_{j0}\approx |I_j|^{-\frac12}\eta_1^5  $ such that
\begin{eqnarray}
&&\|u(t_j)\|_{L^6(|x-x_j|<C(\eta_1) N_j^{-1})}\ge c\eta_1^{\frac
32},\label{25}
\\
&&\|\nabla u(t_j)\|_ {L^2(|x-x_j|<C(\eta_1) N_j^{-1})} \ge
c\eta_1^{\frac 32},\label{26}
\\
&&\|u(t_j)\|_{L^2(|x-x_j|<C(\eta_1) N_j^{-1})}\ge c\eta_1^{\frac
32}N_j^{-1}.\label{27}
\end{eqnarray}

{\bf Proof:} By Beinstein estimate, $\forall N\in 2^Z$, we have
$$
\|P_{\le N} u\|_{\infty} \le N^{\frac12}\|P_{\le N} u\|_6\le
CN^{\frac 12},
$$
which allows us to control the $L^{10}$ norm of low frequency by interpolation,
$$
\|P_{\le N}u\|_{10} \le \|P_{\le N} u\|_{\infty}^{\frac
4{10}}\|P_{\le N}u \|_6^{\frac 6{10}}\le CN^{\frac 15},
$$
hence, using H\"older inequality in time, we have
$$
\|P_{\le N} u\|_{L^{10}(I_j;L^{10})}\le C|I_j|^{\frac 1{10}} N^{\frac 15}.
$$
Taking $N=N_{j0}=C|I_j|^{-\frac 12}\eta_1^5$, one sees that
$$
\|P_{\le N_{j0}} u\|_{L^{10}(I_j;L^{10})}< \frac {\eta_1}2,
$$
and thus
$$
\|P_{\ge N_{j0}} u\|_{L^{10}(I_j;L^{10})}> \frac {\eta_1}2.
$$
Using Littlewood-Paley theorem, we have
\begin{eqnarray*}
(\frac {\eta_1}2)^{10} \le \|P_{\ge N_{j0}}
u\|^{10}_{L^{10}(I_j;L^{10})}&=& \int_{I_j} \|P_{\ge N_{j0}}
u(t)\|_{10}^{10} dt
\\
&=& \int_{I_j} \|(\sum_{N\ge N_{j0}} |P_Nu(t)|^2)^{\frac 12}
\|_{10}^{10}
dt\\
&=&C\int_{I_j}\int_{\mathbf R^3} \sum_{N_1\ge\cdots \ge N_5 \ge
N_{j0}} |P_{N_1} u(t)|^2\cdots |P_{N_5}u(t)|^2 dx dt.
\end{eqnarray*}
Letting $\sigma_N=N^{\frac 12}\|P_N u\|_{L^{\infty}_{xt\in
I_j\times \mathbf R^3}}$, we see the last line is smaller than
\begin{eqnarray}
&&C\sup_{N\ge N_{j0}} \sigma_N^{\frac {20}3}
\int_{I_j}\int_{\mathbf R^3} \sum_{N_1\ge \cdots \ge N_5\ge
N_{j0}}
|P_{N_1}u(t) |^2|P_{N_2} u(t)|^{\frac 43} N_2^{\frac 13} N_3 N_4N_5 dx dt\nonumber\\
&&\le C \sup_{N\ge N_{j0}} \sigma_N^{\frac {20}3}
\int_{I_j}\int_{\mathbf R^3} \sum_{N_1\ge N_2\ge N_{j0}}
N_2^{\frac {10}3} |P_{N_1}u(t)|^2|P_{N_2}u(t)|^{\frac 43} dx dt,
\label{30}
\end{eqnarray}
by summing $N_5$, $N_4$ and $N_3$. Using H\"older inequality and Young's inequality,
(\ref{30}) can be controlled by
\begin{eqnarray*}
&&\quad C\sup_{N\ge N_{j0}} \sigma_N^{\frac {20}3} \sum_{N_1\ge
N_2\ge N_{j0}} N_2^{\frac {10}3} \|P_{N_1}u\|^2_{L_{xt}^{\frac
{10}3}}
\|P_{N_2}  u\|_{L_{xt}^{\frac {10}3}}^{\frac 43} \\
&&\le C\sup_{N\ge N_{j0}} \sigma_N^{\frac {20}3} \sum_{N_1\ge
N_2\ge N_{j0}} N_2^2N_1^{-2}
\|\nabla P_{N_1} u\|^2_{L_{xt}^{\frac {10}3}}\|\nabla P_{N_2} u\|_{L_{xt}^{\frac {10}3}}^{\frac 43}\\
&&\le C\sup_{N\ge N_{j0}}\sigma_{N}^{\frac {20}3} \sum_{N\ge
N_{j0}}\|\nabla P_N
 u\|_{L^{\frac {10}3}(I_j;L^{\frac {10}3})}^{\frac {10}3}\\
&&\le C\sup_{N\ge N_{j0 }} \sigma_N^{\frac {20}3}\|\nabla
u\|^{\frac{10}3}_{L^{\frac {10}3}(I_j;L^{\frac {10}3})}
\\
&&\le C\sup_{N\ge N_{j0}}\sigma_N^{\frac {20}3}.
\end{eqnarray*}
This implies that,
$$
\sup_{N\ge N_{j0}} \sigma_N\ge c\eta_1^{\frac 32},
$$
thus there exists $t_j\in I_j$, $ x_j\in \mathbf R^3$ and $N_j\ge
N_{j0}$ such that
\begin{equation}
|P_{N_j} u(x_j,t_j)|\ge c \eta_1^{\frac 32}N_j^{\frac 12}.
\label{31}
\end{equation}
Now we deduce (\ref{25}), (\ref{26}), (\ref{27}) from (\ref{31}). By the
definition of $P_{N_j}$, we see that
\begin{eqnarray}
c\eta_1^{\frac 32}N_j^{\frac 12} &\le& |P_{N_j}  u(x_j,t_j)| \nonumber \\
&=&|\int \check \phi_{N_j}(x_j-x)  u(t_j,x)dx|\nonumber \\
&\le &\biggl|\int_{|x-x_j|<C(\eta_1) N_j^{-1}} \check \phi_{N_j}(x_j-x)
u(t_j,x)dx\biggr|+
\biggl|\int_{|x-x_j|>C(\eta_1)N_j^{-1}}\check \phi_{N_j}(x_j-x) u(t_j,x)dx\biggr|\nonumber\\
&\le& \biggl(\int_{\mathbf R^3}|\check \phi_{N_j}(x_j-x)|^{\frac
65}dx\biggr)^{\frac 56}
    \biggl(\int_{|x-x_j|<C(\eta_1)N_j^{-1}}| u(t_j,x)|^6dx\biggr)^{\frac 16}\nonumber \\
&\quad& + \biggl(\int_{|x-x_j|>C(\eta_1)N_j^{-1}} |\check
\phi_N(x_j-x)|^{\frac 65}dx\biggr)^{\frac 56}
          \biggl(\int_{\mathbf R^3}| u(t_j,x)|^6dx\biggr)^{\frac 16}.\nonumber \\
\label{32}
\end{eqnarray}
Noting $\check \phi_N(\cdot)=N^3\check\phi(\frac {\cdot}{N})$ and
$\check \phi$ is rapidly decreasing, one obtains
$$
(\ref{32})\le CN_j^{\frac 12}(\int_{|x-x_j|<C(\eta_1)N_j^{-1}}|
u(t,x)|^6)^{\frac 16}+\frac c2 \eta_1^{\frac 32} N_j^{\frac 12},
$$
by choosing $C(\eta_1)$ sufficiently large and
$$
(\int_{\mathbf R^3}| u(t_j,x)|^6dx)^{\frac 16} \le \|
u(t_j)\|_6\le C\Lambda_1^{\frac 16}.
$$
Thus we obtain (\ref{25}). To see (\ref{26}), we begin with
(\ref{31}) that
\begin{eqnarray}
c\eta_1^{\frac 32}N_j^{\frac 12}&<& |P_{N_j} u(x_j,t_j)|=|(\Delta^{-1}\nabla)P_{N_j} \nabla  u(x_j,t_j)|\nonumber\\
&=&|K_{N_j}* \nabla u(x_j,t_j)| =|\int K_{N_j}(x_j-x)\nabla
 u(x,t_j)dx|, \label{33}
\end{eqnarray}
where $K_{N_j}$ is the kernel of $(\Delta^{-1}\nabla)P_{N_j}$,
$\displaystyle K_{N_j}(x)=\mathcal F^{-1}\biggl(\frac {\cdot}{|\cdot|^2}
\phi_{N_j}(\cdot)\biggr) (\xi)$, and
\begin{eqnarray*}
&&\|K_{N_j}\|_{L^2}=N_j^{\frac 12}\|\mathcal F (\frac {\cdot}{|\cdot|^2}\phi(\cdot))\|_2, \\
&&\|K_{N_j}\|_{L^2(|x|\ge C(\eta_1) N_j^{-1})}= N_j^{\frac 12}\|\mathcal
F^{-1}(\frac {\cdot}{|\cdot|^2}\phi)(\cdot)\|_{L^2(|x|\ge
C(\eta_1))} \le C\eta_1^2 N_j^{\frac 12},
\end{eqnarray*}
if $C(\eta_1)$ is
large enough. Thus (\ref{33}) has the bound
\begin{eqnarray*}
&&(\int |K_{N_j}(x_j-x)|^2 dx)^{\frac
12}(\int_{|x-x_j|<C(\eta_1)N_j^{-1}} |\nabla
u(x,t_j)|^2dx)^{\frac 12}
\\
&&\qquad+(\int_{|x-x_j|\ge C(\eta_1)N_j^{-1}}
|K_{N_j}(x_j-x)|^2dx)^{\frac 12}(\int_{\mathbf R^3}|\nabla
u(x,t_j)|^2dx)^
{\frac 12}\\
&&\le CN_j^{\frac 12}(\int_{|x-x_j|<C(\eta_1)N_j^{-1}}|\nabla
u(x,t_j)|^2dx)^{\frac 12} +\frac c2 \eta_1^{\frac 32} N_j^{\frac
12},
\end{eqnarray*}
 and we
have
$$
\|\nabla u(x,t_j)\|_{L^2(|x-x_j|<C(\eta_1)N_j^{-1})}\ge
c\eta_1^{\frac 32}.
$$

The proof of (\ref{27}) is similar. Thus we end the proof of
proposition5.2. Now we use the radial assumption to locate the
bubble at origin.

{\bf Corollary5.3:} Let the conditions in Proposition5.2 be fulfilled. Assume further that $u$ is radial,
then there holds
that
\begin{eqnarray}
&&\|u(t_j)\|_{L^6(|x|<C(\eta_1)N_j^{-1})}\ge c\eta_1^{\frac 32},\label{34}\\
 &&\|\nabla u(t_j)\|_{L^2(|x|<C(\eta_1)N_j^{-1})} \ge c\eta_1^{\frac 32}, \label{35} \\
&&\|u(t_j)\|_{L^2(|x|<C(\eta_1)N_j^{-1})} \ge c\eta_1^{\frac
32}N_j^{-1},\label{36}
\end{eqnarray}
with $t_j$, $N_j$ the same with Proposition5.2.

{\bf Proof:} We prove (\ref{34})--(\ref{36}) by showing that
$$
|x_j|<C(\eta_1)N^{-1}_j,
$$
since once this has been done, we can choose a new constant $\widetilde
C(\eta_1)$ large enough such that
$$B(0,\widetilde C(\eta_1)N^{-1}_j)\supset B(x_j,C(\eta_1) N^{-1}_j),$$
(\ref{34})--(\ref{36}) then
 follow from (\ref{25})--(\ref{27}).

Letting $S(0,|x_j|)$ be a sphere with radius $|x_j|$ and center
$0$. By geometrical observation, one has $O(\frac
{|x_j|}{C(\eta_1) N_j^{-1}})$ balls that have radius $C(\eta_1)
N_j^{-1}$ and center at the points on the sphere. By radial
assumption and Proposition5.2, on each ball, $u(t_j)$ has
nontrivial $L^6$ norm. Using the boundedness of $L^6$ estimate,
one has
$$
O(\frac {|x_j|}{C(\eta_1)N^{-1}}) (C\eta_1^{\frac 32})^6 \le \|u(t)\|_6^6 \le C\Lambda_1.
$$
This gives the desired control on $|x_j|$ and concludes
Corollary5.3.

\section{Proof of Proposition3.5: In case of solitonlike solution}
Applying Corollary5.3 on each interval in the middle component $I^{(2)}$, we get a sequence of time $\{t_j\}$,
$t_j\in I_j$, $\frac {J_1}3+1\le j\le \frac 23 J_1$, such that
\begin{eqnarray}
&&\|\nabla u(t_j)\|_{L^2(|x|\le C(\eta_1)N_j^{-1})} > c\eta_1^{\frac
32},\label{37}\\
&&\|u(t_j) \|_{L^2(|x|\le C(\eta_1)N_j^{-1})} >c\eta_1^{\frac 32}
N_j^{-1}, \quad N_j\ge C|I_j|^{-\frac 12}\eta_1^5.\label{38}
\end{eqnarray}
Now, we discuss two different cases according to the size of the
bubble. First, if there exists $\eta_2$, $0<\eta_2 \ll \eta_1$
such that
\begin{equation}
c|I_j|^{-\frac 12} \eta_1^5 \le N_j\le \frac {C(\eta_1)}{\eta_2} |I_j|^{-\frac 12}, \qquad  \frac {J_1}3+1
\le j\le \frac 23 J_1,
\label{39}
\end{equation}
we call the solution solitionlike. Otherwise there must be $j_0\in [\frac {J_1}3+1,\cdots,\frac 23 J_1]$ such that
\begin{equation}
N_{j_0}\ge \frac {C(\eta_1)}{\eta_2}|I_{j_0}|^{-\frac 12}
\Longleftrightarrow C(\eta_1)N_{j_0}^{-1} <\eta_2|I_{j_0}|^{\frac 12}.
\label{39.1}
\end{equation}
As a consequence, we have concentration as follows,
\begin{equation}
\|\nabla u(t_0)\|_{L^2(|x|<\frac1{\sqrt{2}} \eta_2
|I_{j_0}|^{\frac12})}
> c\eta_1^{\frac 32}. \label{39.2}
\end{equation}
In this case, we call the solution is blow up solution. In this
section, we aim to estimate $J_1$ in case of solitionlike
solution. We follow the idea of \cite{Tao} and begin the proof by
showing that (\ref{38}) holds for every $t\in I_j$, and $\frac
{J_1}3+1\le j\le \frac 23 J_1$.

{\bf Propositin6.1:} Assume $u$ satisfy (\ref{38}), (\ref{39}),
then there exists $C(\eta_1,\eta_2)$, $c(\eta_1,\eta_2)$ such that
\begin{equation}
\|u(t)\|_{L^2(|x|\le C(\eta_1,\eta_2)|I_j|^{\frac 12})}\ge
c(\eta_1,\eta_2)|I_j|^{\frac 12}, \quad \forall t\in I_j \mbox{
and  } j\in [\frac 13 J_1+1,\frac 23 J_1]. \label{40}
\end{equation}
{\bf Proof:} Fix $j$, from (\ref{39}), we have
$$
C\eta_1^{-5}|I_j|^{\frac 12}\ge N_j^{-1}\ge c(\eta_1)\eta_2
|I_j|^{\frac 12}.
$$
Applying this estimate to (\ref{38}), one gets
$$
\|u(t_j)\|_{L^2(|x|<C(\eta_1)|I_j|^{\frac 12})}\ge
c(\eta_1)\eta_2|I_j|^{\frac 12}.
$$
From (\ref{17}) and by choosing $C(\eta_1,\eta_2)$ sufficiently large, we have
\begin{eqnarray*}
\|u(t)\|_{L^2(|x|<C(\eta_1,\eta_2)|I_j|^{\frac 12})}&\ge&
  \|u(t_j)\|_{L^2(|x|<C(\eta_1,\eta_2)|I_j|^{\frac 12})} - \frac {|I_j| \|u\|_{L^{\infty}(I;\dot H^1)}}
  {C(\eta_1,\eta_2)|I_j|^{\frac 12}}\\
&\ge&  c(\eta_1)\eta_2|I_j|^{\frac 12}- c(\eta_1,\eta_2) |I_j|^{\frac 12}\\
&\ge& c(\eta_1,\eta_2)|I_j|^{\frac 12} .
\end{eqnarray*}
This is exactly (\ref{40}).
Once we have gotten (\ref{40}), we can follow the same way in
\cite{Tao} to obtain the finiteness of $J_1$. For the sake of
completeness, we give the proof. First, we do some elementary
computation,
\begin{eqnarray*}
c(\eta_1,\eta_2) |I_j|&\le& \int_{|x|\le C(\eta_1,\eta_2)|I_j|^{\frac 12}}  |u|^2(x,t)dx \\
&\le& \int_{|x|\le C(\eta_1,\eta_2)|I_j|^{\frac 12}} |x|^{\frac 13} \frac {|u|^2(x,t)}{|x|^{\frac 13}}dx \\
&\le& \biggl(\int_{|x|\le C(\eta_1,\eta_2)|I_j|^{\frac 12}} |x|^{\frac
12} dx\biggr)^{\frac 23}
 \biggl(\int_{|x|\le C(\eta_1,\eta_2)|I_j|^{\frac 12}} \frac {|u(x,t)|^6}{|x|}dx\biggr)^{\frac 13}\\
&\le& C(\eta_1,\eta_2)|I_j|^{\frac 76} (\int_{|x|\le
C(\eta_1,\eta_2)|I_j|^{\frac 12}} \frac {|u(x,t)|^6}{|x|} dx)^{\frac
13},
\end{eqnarray*}
thus, we have
\begin{equation}
\int_{|x|\le C(\eta_1,\eta_2)|I_j|^{\frac 12}} \frac {|u(x,t)|^6}{|x|}
dx \ge c(\eta_1,\eta_2) |I_j|^{-\frac 12}. \label{41}
\end{equation}
Comparing (\ref{41}) with Morawetz estimate (\ref{19}), one
obtains,

{\bf Corollary6.2:} For any $I\subset I^{(2)}$, we have
\begin{equation}
\sum_{\frac 13 J_1+1\le j\le \frac 23 J_1; I_j\subset I}
|I_j|^{\frac 12} \le C(\eta_1,\eta_2)|I|^{\frac 12}. \label{42}
\end{equation}
{\bf Proof:} Noting $|I_j|^{\frac 12}<|I|^{\frac 12}$ and letting
$A=C(\eta_1,\eta_2)$, (\ref{41}) becomes
\begin{equation}
\int_{|x|\le A|I|^{\frac 12}} \frac {|u|^6(x,t)}{|x|} dx \ge
c(\eta_1,\eta_2)|I_j|^{-\frac 12}. \label{43}
\end{equation}
Integrating (\ref{43}) on $I_j$ and summing together in $j$, we
get,
\begin{eqnarray*}
c(\eta_1,\eta_2)\sum_{\frac 13 J_1+1\le j\le \frac 23 J_1; I_j\subset
I} |I_j|^{\frac 12}
&\le& \int_I \int_{|x|\le A|I|^{\frac 12}}\frac {|u|^6(x,t)}{|x|} dx dt \\
&\le & CA|I|^{\frac 12}\\
&\le & C(\eta_1,\eta_2) |I|^{\frac 12},
\end{eqnarray*}
this gives (\ref{42}). As a direct consequence of Corollary6.2, we
have

{\bf Corollary6.3:} Let $\displaystyle I=\cup_{j_1\le j\le j_2}
I_j$ be a union of consecutive intervals, $\frac 13 J_1+1\le
j_1,j_2\le \frac 23 J_1$, then there exists $j_1\le j\le j_2$ such
that $|I_j|>c(\eta_1,\eta_2)|I|$.

 {\bf Proof:} From (\ref{42}) we
know that
$$
C(\eta_1,\eta_2)|I|^{\frac 12} \ge \sum_{j_1\le j\le j_2}
|I_j|^{\frac 12}\ge \sum_{j_1\le j\le j_2} |I_j| (\sup_{j_1\le
j\le j_2} |I_j|)^{-\frac 12}=|I|(\sup_{j_1\le j\le
j_2}|I_j|)^{-\frac 12} ,
$$
and hence
\begin{equation}
C(\eta_1,\eta_2)|I|^{-\frac 12}\ge (\sup_{j_1\le j\le j_2}
|I_j|)^{-\frac 12}. \label{44}
\end{equation}
(\ref{44}) allows us to find an interval $I_j$ such that
$$
|I_j|>c(\eta_1,\eta_2)|I| .
$$
Now, we show that the intervals $I_j$ must concentrate at some
time $t_*$.

{\bf Proposition6.4:} There exists $t_*\in I^{(2)}$ and distinct
intervals $I_{j_1},\cdots,I_{j_K}$, $j_k \in [\frac
{J_1}3+1,\cdots\frac 23 J_1]$, $K>C(\eta_1,\eta_2)$ such that
$$
|I_{j_1}|\ge 2|I_{j_2}|\ge \cdots \ge 2^{k-1}|I_{j_K}|,$$ and
$dist(t_*,I_{j_k}) \le C(\eta_1,\eta_2)|I_{j_k}|$.

For the proof of this Proposition, one refers to \cite{Tao},Proposition3.8.

Let $t_*$ and $I_{j_1},\cdots, I_{j_k}, \cdots,
I_{j_K}$ be as in the Proposition6.4 and for every $t\in I_{j_k}$,
there holds
\begin{equation}
Mass(u(t),B(0; C(\eta_1,\eta_2)|I_{j_k}|^{\frac 12})) \ge
c(\eta_1,\eta_2) |I_{j_k}|^{\frac 12},\quad \forall t\ \in
I_{j_k}, \qquad 1\le k\le K.\label{g}
\end{equation}
The point $0$ can be substituted by $x_{j_k}$ without modification
to the following proof, as just like the thing that has been
mentioned in \cite{Tao}. By the local mass conservation, we have
\begin{eqnarray*}
Mass(u(t_*),B(0; C(\eta_1,\eta_2)))|I_{j_k}|^{\frac 12} &\ge&
c(\eta_1,\eta_2)|I_{j_k}|^{\frac 12}-\frac {|t_*-t|
\|u\|_{L^{\infty}(I_{j_k},\dot H^1)}}
{C(\eta_1,\eta_2)|I_{j_k}|^{\frac 12}},\\
&\ge& c(\eta_1,\eta_2) |I_{j_k}|^{\frac 12}, \qquad \forall 1\le
k\le K.
\end{eqnarray*}
Denote $B_k=B(0; C(\eta_1,\eta_2) |I_{j_k}|^{\frac 12})$, we
rewrite the above estimate as follows,
\begin{equation}
Mass(u(t_*),B_k)\ge c(\eta_1,\eta_2)|I_{j_k}|^{\frac 12}, \qquad
1\le k\le K .\label{45}
\end{equation}
On the other hand, by the local mass estimate (\ref{18}), we have
$$
Mass(u(t_*),B_k)\le C(\eta_1,\eta_2)|I_{j_k}|^{\frac 12}.
$$
Letting $N:= \log(\frac 1{\eta_3})$, then for $k'>k+N$, we have that
\begin{eqnarray*}
\int_{B'_k} |u(t_*,x)|^2 dx \le C(\eta_1,\eta_2)|I_{j_k'}|&\le& C(\eta_1,\eta_2) 2^{-(k'-k)}|I_{j_k}|\\
&\le& C(\eta_1,\eta_2) \eta_3 2^{-(k'-k-N)} |I_{j_k}|,
\end{eqnarray*}
and hence,
\begin{equation}
\sum_{k+N\le k'\le K} \int_{B_{k'}}|u(t_*,x)|^2dx\le
C(\eta_1,\eta_2)\eta_3 |I_{j_k}| \sum_{k_N \le k'\le  K}
2^{-(k'-k-N)}. \label{46}
\end{equation}
By the finiteness of the summation, the assumption on
$\eta_1,\eta_2,\eta_3$, and (\ref{g}), we continue to estimate
(\ref{46}) by
$$
c(\eta_1,\eta_2)|I_{j_k}|\le \frac 12 Mass (u(t_*),B_k)^2=\frac
12\int_{B_k} |u(t_*,x)|^2dx,
$$
and hence
\begin{eqnarray}
\int_{B_k \backslash (\cup_{k+N\le k'\le  K} B_{k'})} |u(t_*,x)|^2dx &\ge&
\int_{B_k} |u(t_*,x)|^2 dx -\int_{\cup_{k+N\le k'\le K} B_{k'}}
|u(t_*,x)|^2 dx \nonumber\\
&\ge& \int_{B_k} |u(t_*,x)|^2dx -\sum_{k+N\le k' \le K} \int_{B_{k'}} |u(t_*,x)|^2 dx \nonumber \\
&\ge& \frac 12 \int_{B_k}|u(t_*,x)|^2 dx \ge
c(\eta_1,\eta_2)|I_{j_k}|. \label{46.1}
\end{eqnarray}
By H\"older inequality, we further give the upper bounds of the
left side as follows,
\begin{eqnarray*}
&&\qquad \int_{B_k\backslash (\cup_{k+N\le k'\le K} B'_k)} |u(t_*,x)|^2dx\\
&\le&  (\int_{B_k\backslash (\cup_{k+N\le k'\le K} B'_k)} |u(t_*,x)|^6 dx)^{\frac 13} mes(B_k)^{\frac 23} \\
&\le& C(\eta_1,\eta_2)|I_{j_k}| (\int_{B_k\backslash (\cup_{k+N\le
k'\le K} B'_k)} |u(t_*,x)|^6dx)^{\frac 13},
\end{eqnarray*}
hence we have
\begin{equation}
\int_{B_k\backslash (\cup_{k+N\le k'\le K} B'_k)}
|u(t_*,x)|^6dx \ge c(\eta_1,\eta_2). \label{47}
\end{equation}
Summing (\ref{47}) in $k$, we obtain
\begin{equation}
\sum_{k=1}^{ K} \int_{B_k\backslash (\cup_{k+N\le k'\le
K} B'_k)} |u(t_*,x)|^6dx \ge c(\eta_1,\eta_2) K. \label{48}
\end{equation}
Denoting $P_k:= B_k\backslash (\cup_{k+N\le k'\le  K} B'_k)$, then
$\{ P_k\}_{k=1}^{ K}$ overlaps at most $N$ times. Thus the
left hand side of (\ref{48}) is smaller than
$$
N\int_{\mathbf R^3} |u(t_*,x)|^6dx.
$$
By the definition of $\eta_3$, the boundedness of $\|u(t_*)\|_6$,
we have an upper bound for $  K$,
$$
 K\le C(\eta_1,\eta_2,\eta_3,\Lambda_1,\Lambda_2),
$$
and this in turn gives the control of $J_1$,
$$
J_1\le C \exp (C(\eta_1,\eta_2,\eta_3,\Lambda_1,\Lambda_2)).
$$

\section{In case of blow up solution}
Our purpose of this section is to prove the boundedness of $J_1$
under the condition (\ref{39.1}) and (\ref{39.2}). That is ,for
the solution, we have  concentration at some $t_0\in I_{j_0}$,
$j_0\in[\frac {J_1}3+1,\cdots,\frac 23 J_2]$ such that
\begin{equation}
\|\nabla u(t_0)\|_{L^2(|x|<\frac 1{\sqrt{2}} \eta_2 |I_{j_0}|^{\frac12})}>
c\eta_1^{\frac 32} .\label{50}
\end{equation}
If $t_0$ lies in the left side of $I_{j_0}$, we take $I=[t_0,b]$
where $b$ is the left end point of $I_{j_0}$;Otherwise we take
$I=[a,t_0]$ with $a$ the right end point of $I_{j_0}$. Then
(\ref{50}) becomes
\begin{equation}
\|\nabla u(t_0)\|_{L^2(|x|<\eta_2 |I|^{\frac 12})}>c\eta_1^{\frac
32}. \label{51}
\end{equation}
Assume $I=[t_0,b]$, we plan to re-solve the problem (\ref{1.1})
forward in time. Otherwise, we do in the reverse direction. First
we show that, by removing the small bubble, we remove nontrivial
portion of energy.

Let $\chi$ be a smooth radial function such that $\chi(x)=1$ as $|x|\le1$, and $\chi(x)=0$ as
$|x|\ge 2$. Let $\phi(x)=\chi\biggl(\displaystyle{\frac{x}{N\eta_2|I|^{\frac12}}}\biggr)$ for some $N\ge1$ to
be specified later, and $w(t_0,x)=(1-\phi(x))u(t_0,x)$, then we have


\bf{Lemma7.1} \rm $E_1(w(t_0))\le E_1(u(t_0))-c\eta_1^3$.

\bf{Proof:} \rm By noting
$$
w(t_0)=(1-\phi)u(t_0),
$$
we compute that
$$
\nabla w(t_0)=(1-\phi)\nabla u(t_0)-\nabla\phi u(t_0),
$$
and thus,
$$
|\nabla w(t_0)|^2=|\nabla u(t_0)|^2+(\phi^2-2\phi)|\nabla u(t_0)|^2+
|\nabla\phi|^2|u(t_0)|^2-2 Re (1-\phi)\nabla\phi\bar u(t_0)\nabla u(t_0).
$$
Integrating it on $\mathbf R^3$, one gets
\begin{eqnarray*}
\|\nabla w(t_0)\|_2^2&\le& \|\nabla u(t_0)\|_2^2+\int_{\mathbf R^3}(\phi^2-2\phi)|\nabla
u(t_0,x)|^2dx\\
&&\quad -2\int_{\mathbf R^3}|\nabla \phi(x)
u(t_0,x)|^2dx-2Re\int_{\mathbf R^3}(1-\phi)\nabla\phi \bar
u(t_0)\nabla u(t_0)(x)dx.
\end{eqnarray*}
By the trivial inequality: $\phi^2-2\phi\le-\phi$ and (\ref{51}),
one can estimate the second term of the right side by
$$
-\int _{|x|\le N\eta_2|I|^{\frac12}}|\nabla u(t_0,x)|^2 dx\le
-c\eta_1^3.
$$
Now, we estimate the remaining two terms. We use H\"older inequality to control them by
\begin{eqnarray}
&&C\|\nabla \phi\|_3^2\|u(t_0)\|_{L^6_{N\eta_2|I|^{\frac12}\le
|x|\le 2N\eta_2|I|^{\frac12}}}^2+C\|\nabla\phi\|_3\|\nabla
u(t_0)\|_2 \|u(t_0)\|_{L^6_{N\eta_2|I|^{\frac12}\le|x|\le
2N\eta_2|I|^{\frac12}}}\nonumber\\
&\le& C\biggl(\|u(t_0)\|_{L^6_{N\eta_2|I|^{\frac12}\le |x|\le 2N\eta_2|I|^{\frac12}}}^2+\|u(t_0)\|_
{L^6_{N\eta_2|I|^{\frac12}\le |x|\le 2N\eta_2|I|^{\frac12}}}\|\nabla
u(t_0)\|_2\biggr).\nonumber\\
\label{56}
\end{eqnarray}
Now, we claim that, there must exist $N$ which depend only on
$\eta_1$ such that
\begin{equation}
\|u(t_0)\|_{L^6_{N\eta_2|I|^{\frac12}\le|x|\le
2N\eta_2|I|^{\frac12}}}\le \eta_1^4.\label{56.1}
\end{equation}
Indeed, if otherwise, we will have $N$ annuluses , on each annulus, $u(t_0)$ has nontrivial $L^6$ norm. Summing
these annuluses together, we obtain
$$
N(\eta_1^4)^6\le \sum_{N'\le N\in \mathbf
N}\|u(t_0)\|_{L^6_{N'\eta_2|I|^{\frac12}\le|x|\le
2N'\eta_2|I|^{\frac12}}}\le C,
$$
by the boundedness of $L^6$ estimate. This will be a contradiction if $N\ge C\eta_1^{-24}$. Hence, one can fix $N=
C(\eta_1)$ such that (\ref{56.1}) holds and
$$
(\ref{56})\le C\eta_1^4.
$$
We finally obtain this Lemma by noting
$$
E_1(w(t_0))=\frac 12\|\nabla w(t_0)\|_2^2+\frac 13\|w(t_0)\|_6^6,
$$
and combing the above estimates together.

\bf{Lemma7.2} \rm  We have that,
\begin{eqnarray*}
&&E_1(w(t_0))\le \Lambda_1-c\eta_1^3,\\
&&E_2(w(t_0))\le \Lambda_2+C\eta_1^4.
\end{eqnarray*}
\bf{Proof:} \rm Noting Lemma7.1, it suffices to prove
\begin{eqnarray*}
&&E_1(u(t_0))\le E_1(u(0))+C\eta_1^4,\\
 &&E_2(u(t_0))\le E_2(u(0))+C\eta_1^4.
 \end{eqnarray*}
So, Let's compute the increment of $E_i(u(t))$ from 0 to $t_0$:
$$
\int_0^{t_0}\frac{\partial}{\partial t}E_1(u(t))dt, \mbox{ and }\int_0^{t_0}\frac{\partial}{\partial t}E_2(u(t))dt.
$$
From the equation (\ref{1.1}), we see that
\begin{eqnarray*}
\frac{\partial}{\partial t} E_2(u(t))&=&\frac{\partial}{\partial
t}\|xu(t)\|_2^2\\
&=& 2Im\int_{\mathbf R^3} x\bar u\nabla u(t,x)dx\\
&\le& C\|xu\|_{L^{\infty}((0,t_0);L^2)}\|\nabla u\|_{L^{\infty}((0,t_0);L^2)}\le C(\Lambda_1,\Lambda_2).\quad \forall
t\in [0,t_0],
\end{eqnarray*}
thus, we have
$$
\biggl|\int_0^{t_0}\frac{\partial}{\partial t}E_2(t)dt\biggr|\le C\eta_1^4.
$$
By noting that $\frac{\partial}{\partial t}E_1(t)=-\frac{\partial}{\partial t}E_2(t)$, we get
$$
\biggl|\int_0^{t_0}\frac{\partial}{\partial t}E_1(t)dt\biggr|\le C\eta_1^4,
$$
hence, Lemma7.2 follows.

Putting these Lemmas aside, we turn to re-solve the solution from $t_0$ forward. We do this by splitting
$u=v+w$ and studying the following two initial data problems:
\begin{equation}
\left\{
\begin{array}{l}
(i\partial_t+\frac{\Delta}2+\frac{|x|^2}2)v=|v|^4v,\\
v(x,t_0)=\phi(x)u(t_0,x).
\end{array}
\right.
\label{57}
\end{equation}
\begin{equation}
\left\{
\begin{array}{l}
(i\partial_t+\frac{\Delta}2+\frac{|x|^2}2)w=|v+w|^4(v+w)-|v|^4v,\\
w(x,t_0)=(1-\phi(x))u(t_0,x).
\end{array}
\right.
\label{58}
\end{equation}
Now, Let's first prove that (\ref{57}) is  wellposed on $[t_0,\infty)$.

\bf{Proposition7.3} \rm There exists a unique solution $v(x,t)$ to (\ref{57}) satisfies
\begin{eqnarray*}
&&\|v\|_{L^{10}(I;L^{10})}\le C\eta_1, \quad
\|v\|_{L^{10}([b,\infty);L^{10})}\le C\eta_2^{\frac
15}.\\
&&\|A(\cdot)v\|_{L^q(I;L^r)} \le C,
\end{eqnarray*}
where $A\in\{J,H\}$ and $(q,r)$ are admissible pairs.

\bf{Proof:} \rm We begin by computing the $L^{10}$ norm of the linear flow $U(t-t_0)(\phi u(t_0))$. First, by Duhamel's
formula, we observe that
$$
U(t-t_0)u(t_0)=u(t)+i\int_{t_0}^t U(t-s)|u|^4u(s)ds,
$$
from this, we see that
$$
\|U(\cdot-t_0)u(t_0)\|_{L^{10}(I;L^{10})}\le \|u\|_{L^{10}(I;L^{10})}
+\|\int_{t_0}^t U(t-s)|u|^4u(s)ds\|_{L^{10}(I;L^{10})}.
$$
Applying embedding and Strichartz, the second term is smaller than
\begin{eqnarray*}
&&\|J(\cdot)\int_{t_0}^t
U(t-s)|u|^4u(s)ds\|_{L^{10}(I;L^{\frac{30}{13}})}\\
&\le&
C\|u\|_{L^{10}(I;L^{10})}^4\|J(\cdot)u\|_{L^{\frac{10}3}(I;L^{\frac{10}3})}\\
&\le& C\eta_1^4\le \eta_1,
\end{eqnarray*}
and hence,
$$
\|U(\cdot-t_0)u(t_0)\|_{L^{10}(I;L^{10})}\le 2\eta_1.
$$
Noting $\phi u(t_0)$ is a radial function in space, we have that
\begin{eqnarray*}
U(t-t_0)(\phi u(t_0))(x)&=&\exp{\{\frac{it(\Delta+|x|^2)}2\}}(\phi
u(t_0))(x)\\
&=&\mathcal F^{-1}\biggl(\exp{\{-\frac{it(|\xi|^2+\Delta_{\xi})}2\}}\widehat{\phi
u(t_0)}(\xi)\biggr)(x)\\
&=&\int_{\mathbf R^3} e^{ix\xi}e^{-\frac{i(t-t_0)}2(|\xi|^2+\Delta_{\xi})}\int_{\mathbf R^3}\hat u(t_0)(\xi-\xi_1)\hat\phi(\xi_1)d\xi_1d\xi.
\end{eqnarray*}
Expanding $|\xi|^2=|\xi-\xi_1|^2+2\xi_1(\xi-\xi_1)+|\xi_1|^2$, the above term becomes
$$
\int\int e^{i(\xi-\xi_1)(x-(t-t_0)\xi_1)}e^{-\frac{i(t-t_0)}2(|\xi-\xi_1|^2+\Delta_{\xi-\xi_1})}\hat u(t_0)(\xi-\xi_1)d\xi
e^{-\frac{i(t-t_0)}2|\xi_1|^2 }\hat \phi(\xi_1)d\xi_1 ,
$$
by renaming the variable, one sees that this is exactly
$$
\int_{\mathbf R^3} U(t-t_0)u(t_0)(x-(t-t_0)\xi_1)e^{-\frac{i(t-t_0)}2|\xi|^2}\hat \phi(\xi_1)d\xi_1,
$$
and hence
\begin{eqnarray*}
\|U(\cdot-t_0)(\phi u(t_0))\|_{L^{10}(I;L^{10})}&\le &\|U(\cdot-t_0)u(t_0)\|_{L^{10}(I;L^{10})}\|\hat
\phi\|_1\\
&\le& C\|U(\cdot-t_0)u(t_0)\|_{L^{10}(I;L^{10})}
\le C\eta_1.
\end{eqnarray*}
The estimate of the linear flow allows us to solve the problem in the following set,
$$
X:=\biggl\{v(x,t)\biggl|\quad \|v\|_{L^{10}(I;L^{10})}\le
C\eta_1,\quad \|A(\cdot)v\|_{L^{\frac{10}3}(I;L^{\frac{10}3})}\le
C,\quad A\in\{J,H\} \biggr\},
$$
donated with the metric
$$
d(u_1,u_2)=\|u_1-u_2\|_{L^{10}(I;L^{10})}+\max_{A\in\{J,H\}}\|A(\cdot)(u_1-u_2)\|_{L^{\frac{10}3}(I;L^{\frac{10}3})}.
$$
We omit the proof of this part since it is routine. Once we have
gotten the solution on $I=[t_0,b]$, we extend this solution beyond
$I$. As a consequence, we are left to show a finite apriori
spacetime estimate on $[b,\infty)$. Assuming $v$ be a finite
energy solution on $[b,\infty)$, we redefine the energy of $v$ by
\begin{eqnarray*}
\tilde\mathcal E_1(t)=\frac12\|J(t-t_0)v(t)\|_2^2+\frac
13\cosh^2(t-t_0)\|v(t)\|_6^6;\\
\tilde\mathcal E_1(t)=\frac12\|H(t-t_0)v(t)\|_2^2+\frac
13\sinh^2(t-t_0)\|v(t)\|_6^6;
\end{eqnarray*}
Repeating the computations in Lemma2.4, we find
$$
\frac{d\tilde \mathcal E_1(t)}{dt}=-\frac23\sinh2(t-t_0)\|v(t)\|_6^6=\frac{d\tilde \mathcal E_2(t)}{dt}.
$$
Integrating the second half of the equation, we have
\begin{eqnarray*}
&&\frac12\|H(t-t_0)v(t)\|_2^2+\frac13\sinh^2(t-t_0)\|v(t)\|_6^6\\
&=& \frac12 \|xv(t_0)\|_2^2-\frac23\int_{t_0}^t\sinh(2(\tau-t_0))\|v(\tau)\|_6^6d\tau
\end{eqnarray*}
This implies
\begin{equation}
\sinh^2(t-t_0)\|v(t)\|_6^6\le C\|x\phi u(t_0)\|_2^2.
\label{59}
\end{equation}
By H\"older and direct computation, we continue to estimate the right side as
\begin{eqnarray*}
\|x\phi u(t_0)\|_2&\le&
\|\cdot\chi(\frac{\cdot}{C(\eta_1)\eta_2|I|^{\frac12}})u(t_0)\|_2\\
&\le& (C(\eta_1)\eta_2)^2|I|\|\cdot\chi(\cdot)\|_3\|u(t_0)\|_6\\
&\le& C(\eta_1)\eta_2^2|I|.\\
\sinh^2(t-t_0)  &\ge&|t-t_0|^2\le |I|^2, \quad t\ge b,
\end{eqnarray*}
hence, from (\ref{59}), we have
\begin{equation}
\|v(t)\|_6^6\le C(\eta_1)\eta_2^4\le \eta_2^2,\quad\forall t\ge b.
\label{60}
\end{equation}
On the other hand, noting on $[b,\infty)$, $v$ satisfies,
\begin{equation}
v(t)=U(t-t_0)v(t_0)-i\int_{t_0}^tU(t-s)|v|^4v(s)ds,
\label{61}
\end{equation}
we have
\begin{eqnarray*}
\|J(\cdot)v\|_{L^6([b,\infty);L^{\frac{18}7})}&\le& C\|J(t_0)v(t_0)\|_2+C\|J(\cdot)|v|^4v\|
_{L^{\frac32}([b,\infty);L^{\frac{18}{13}})}\\
&\le & C\|J(t_0)v(t_0)\|_2+C\|v\|_{L^{\infty}([b,\infty);L^{6})}\|v\|^3_{L^6([b,\infty);L^{18})}
 \|J(\cdot)v\|_{L^6([b,\infty);L^{\frac{18}7})}\\
 &\le&
 C\|J(t_0)v(t_0)\|_2+C\eta_2^{\frac13}\|J(\cdot)v\|^4_{L^6([b,\infty);L^{\frac{18}7})}.
 \end{eqnarray*}
This implies
\begin{equation}
\|J(\cdot)v\|_{L^6([b,\infty);L^{\frac{18}7})}\le C\|J(t_0)v(t_0)\|_2\le C,
\label{62}
\end{equation}
where, $C$ depends only on $\Lambda_1$, $\Lambda_2$. To see this, we use Lemma2.3 to expand $J(t_0)v(t_0)$ as
$$
\phi(x)J(t_0)u(t_0)+i\cosh t_0u(t_0)\nabla \phi(x),
$$
which can be easily controlled.

Combing the bounds (\ref{60}) and (\ref{62}) together and using interpolation, one obtains
\begin{eqnarray}
\|v\|_{L^{10}([b,\infty);L^{10})}&\le& C\|v\|^{\frac4{10}}_{L^{\infty}([b,\infty);L^6)}\|J(\cdot)v\|^{\frac 6{10}}
_{L^6([b,\infty);L^{\frac{18}7})}\nonumber\\
&\le & C\eta_2^{\frac15}.\label{62.1}
\end{eqnarray}
This combining with some routine arguments gives Proposition7.3.

Now, we are at the position to solve the Cauchy problem (\ref{58}). Before doing this, we list the estimates
that follows from Proposition7.3 and the conditions on $u$.
\begin{equation}
\begin{array}{c}
\|v\|_{L^{10}([b,\infty);L^{10})}<C\eta_2^{\frac15}\\
\|w\|_{L^{10}(I;L^{10})}\le C\eta_1,\quad \|A(\cdot)w\|_{L^{\frac{10}3}(I;L^{\frac{10}3})}\le C;\quad
 \quad \|B w\|_{L^{\frac{10}3}(I;L^{\frac{10}3})}\le C,\\
\|v\|_{L^{10}(I;L^{10})}\le C\eta_1,\quad \|A(\cdot)v\|_{L^{\frac{10}3}(I;L^{\frac{10}3})}\le C;\quad
 \quad \|B v\|_{L^{\frac{10}3}(I;L^{\frac{10}3})}\le C,\\
 \qquad\qquad A\in\{J,H\},\quad B\in\{i\nabla_x,x\},
\end{array}
\label{63}
\end{equation}
here, we have used the condition that $I\in[0, \eta_1^4)$ to get the estimate on $Bv$, $Bw$. The constants
above depend only on $\Lambda_1$, $\Lambda_2$.

For the sake of doing perturbation analysis and applying the induction, it's necessary to introduce the following Lemma.

\bf{Lemma7.4}. \rm We have that
$$E_1(w(b))\le \Lambda_1-c\eta_1^3,\quad E_2(w(b))\le \Lambda_2+C\eta_1^4.$$
\bf{Proof:} \rm Noting  Lemma7.2, we need only to prove that
\begin{eqnarray*}
&&\biggl|\int_{t_0}^b\frac{\partial}{\partial t} E_1(w(t))dt\biggr|\le
C\eta_1^4,\\
&&\biggl|\int_{t_0}^b\frac{\partial}{\partial t} E_2(w(t))dt\biggr|\le
C\eta_1^4,
\end{eqnarray*}
For simplicity, denote
$$
|v+w|^4(v+w)-|v|^4v=|w|^4w+F(v,w),
$$
hence, $w$ satisfies the equation
$$
(i\partial_t+\frac{\Delta}2+\frac{|x|^2}2)w=|w|^4w+F(v,w).
$$
By some basic computation, one sees that
$$
\frac{\partial}{\partial_t} E_2(w(t))= 2Im\int_{\mathbf R^3} x\bar w\nabla w dx+2Im\int_{\mathbf R^3}
|x|^2\bar wF(v,w)dx,
$$
thus, we get
\begin{eqnarray*}
\biggl|\int_{t_0}^b \frac{\partial}{\partial t} E_2(w(t))dt\biggr|&\le& 2|t_0-b|\|xw\|_{L^{\infty}(I;L^2)}\|
\nabla w \|_{L^{\infty}(I;L^2)}\\
&&+C\|xw\|^2_{L^{\frac{10}3}(I;L^{\frac{10}3})}(\|w\|_{L^{10}(I;L^{10})}^4+\|v\|^4_{L^{10}(I;L^{10})})\le C\eta_1^4.
\end{eqnarray*}
To prove the increment of the $E_1(w(t))$ from $t_0$ to $b$, we first compute directly that
\begin{eqnarray*}
\frac{\partial}{\partial t} E_1(w(t))&=&\frac{\partial}{\partial t}E_2(w(t))+Re \int_{\mathbf R^3} F(v,w)\bar w_t(x)
dx\\
&=&\frac{\partial}{\partial t} E_2(w(t))+Im\int_{\mathbf R^3}
F(v,w)(\frac12\Delta\bar w+\frac12|x|^2\bar w-|w|^4\bar w-\overline{F(v,w)})(x)dx
\end{eqnarray*}
Integrating over $[t_0,b]$ and using integration by parts, one gets
\begin{eqnarray*}
\biggl|\int_{t_0}^b \frac{\partial}{\partial t} E_1(w(t))\biggr|&\le& C\eta_1^4
+ C(\|\nabla
w\|_{L^{\frac{10}3}(I;L^{\frac{10}3})}^2+\|xw\|^2_{L^{\frac{10}3}(I;L^{\frac{10}3})})\\
&&\times (\|v\|_{L^{10}(I;L^{10})}^4+\|w\|_{L^{10}(I;L^{10})}^4)+C(\|v\|_{L^{10}(I;L^{10})}^{10}+
\|w\|_{L^{10}(I;L^{10})}^{10})\\
&\le & C\eta_1^4.
\end{eqnarray*}
This ends Lemma7.4.

Now, for the sake of convenience, we make a small adjustment such that,
the increment of $E_1$ and the decrement $E_2$ take the same value. More
precisely, noting Lemma7.4, we can get
$$
E_1(w(b))\le \Lambda_1-C\eta_1^4,\quad E_2(w(b))\le \Lambda_2+C\eta_1^4,
$$
here, the above two constants are same.

Now, we make an induction assumption in order to solve the problem (\ref{58}). We assume that:

Let $t'\in \mathbf R$ and $W(t',x)$ satisfy
$$
E_1(W(t'))\le \Lambda_1-C \eta_1^4, \quad E_2(W(t'))\le \Lambda_2+C \eta_1^4.
$$
Then the Cauchy problem of (\ref{1.1}) with prescribed data $W(t')$ at $t'$ is solvable on $[t'-\eta_1^4,t'+\eta_1^4]$,
and the solution $W$ satisfies
$$
\|W\|_{L^{10}([t'-\eta_1^4 ,t'+\eta_1^4];L^{10})}\le C(\Lambda_1-C\eta_1^4,\Lambda_2+C\eta_1^4).
$$

By this assumption and Lemma7.4, we see that the solution of
\begin{eqnarray*}
\left\{
\begin{array}{ll}
i W_t+\frac {\Delta }2 W=-\frac {|x|^2}2 W+|W|^4W,\\
W(b)=w(b)
\end{array}
\right.
\end{eqnarray*}
satisfies the estimate
$$
\|W\|_{L^{10}([b-\eta_1^4,b+\eta_1^4];L^{10})}\le C(\Lambda_1-C\eta_1^4,\Lambda_2+C\eta_1^4)\le C(\Lambda_1,\Lambda_2).
$$
Substracting $W$ from $w$, we are left to solve the perturbation problem with respect to $\Gamma=w-W$ on $[b,\eta_1^4]$,
\begin{equation}
\left\{
\begin{array}{ll}
(i\partial_t+\frac {\Delta}2+\frac
{|x|^2}2)\Gamma=|v+W+\Gamma|^4(v+W+\Gamma)-|v|^4v-|W|^4W,\\
\Gamma(b)=0 .
\end{array}
\right.
\label{70}
\end{equation}

\section{Solving the perturbation problem}
Our task of this section is to solve (\ref{70}) with the help of (\ref{63}). To insure the smallness of the nonlinear flow,
we split $[b,\eta_1^4]$ into finite subintervals such that on each subinterval, $W$ is small, so that we can solve (\ref{70}) on
every subinterval. Before doing this, we re-estimate $v$ on $[b,\infty)$.

{\bf Lemma8.1:} Excepting for (\ref{63}), $v$ satisfies
$$
\|A(\cdot) v\|_{L^{\frac {10}3}([b,\infty);L^{\frac {10}3})}\le c(\eta_2),\qquad A\in\{J,H\}.
$$

{\bf Proof:} Taking $J$ as an example, one sees
$$
J(t)v(t)=U(t-t_0) J(t_0)v(t_0)-i\int_{t_0}^t U(t-s)J(s)|v|^4 v(s)ds.
$$
For the linear term, we estimate directly. From decay estimate (\ref{7.1}),
\begin{eqnarray*}
\|U(t-t_0)J(t_0)v(t_0)\|_{L^{\infty}_x} &\le&
C|t-t_0|^{-\frac32}\|J(t_0)v(t_0)\|_1,\\
&\le& C|t-t_0|^{-\frac 32}(\|J(t_0)u(t_0)\phi\|_1+\|\cosh t_0 u(t_0) \nabla_x \phi \|_1)
\\
&\le& C|t-t_0|^{-\frac 32}(\|\phi\|_2\|J(t_0)u(t_0)\|_2+\|\cosh t_0
u(t_0)\|_6\|\nabla_x \phi\|_{\frac 65})
\end{eqnarray*}
By noting $\phi(x)=\chi(\frac x{\eta_1 C(\eta_1)|I|^{\frac 12}})$, one has
\begin{equation}
\|U(t-t_0)J(t_0) v(t_0)\|_{L^{\infty}_x} \le C(\eta_1)(\frac {\eta_2|I|^{\frac 12}}{t-t_0})^{\frac 32}.
\label{71}
\end{equation}
On the other hand,
\begin{eqnarray}
\|U(t-t_0)J(t_0) v(t_0)\|_{L^2_x}&\le& \|J(t_0)v(t_0)\|_{L^2}\nonumber \\
&\le& \|J(t_0)u(t_0)\|_2\|\phi\|_{\infty}+\|\nabla \phi \|_3 \|\cosh t_0
u(t_0)\|_6 \nonumber\\
&\le& C. \label{72}
\end{eqnarray}
By interpolation and (\ref{71}) and (\ref{72}), we have that
\begin{eqnarray*}
\|U(t-t_0)J(t_0)v(t_0)\|_{L^{\frac {10}3}_x} &\le&
\|U(t-t_0)J(t_0) v(t_0)\|_{\infty}^{\frac
25}\|U(t-t_0)J(t_0)v(t_0)\|_2^{\frac 35}
\\
&\le& C(\eta_1)(\frac {|I|^{\frac 12} \eta_2}{|t-t_0|})^{\frac 35},
\end{eqnarray*}
and thus
\begin{eqnarray}
\|U(t-t_0)J(t_0) v(t_0)\|^{\frac {10}3}_{L^{\frac {10}3}([b,\infty);L^{\frac {10}3})}&\le&
C(\eta_1)|I| \eta_2^2 \int_b^{\infty} \frac {dt}{|t-t_0|^2}\nonumber \\
&\le& C(\eta_1)\eta_2^2 .\label{73}
\end{eqnarray}
To estimate the nonlinear term, we denote $t_1=t_0+\eta_2|I|$, and split it into two parts,
$$
\|\int_{t_0}^{t_1}U(t-s) J(s) |u|^4 u(s) ds\|_{L^{\frac {10}3}([b,\infty);L^{\frac {10}3})}+
\|\int_{t_1}^{\infty}U(t-s) J(s) |u|^4 u(s) ds\|_{L^{\frac {10}3}([b,\infty);L^{\frac {10}3})},
$$
For the first part, we use $L^p - L^{p'}$ estimate to control it by
$$
\|\int_{t_0}^{t_1} |t-s|^{-\frac 35}\|J(s)|u|^4 u(s)\|_{L^{\frac {10}7}} ds \|_{L^{\frac {10}3}([b,\infty)}.
$$
Since for $s\in [t_0,t_1]$, $t>b$, $|t-s|\sim |t-t_0|$, we see the first part is smaller than
\begin{eqnarray*}
&&\quad C\||t-t_0|^{-\frac 53}\|_{L^{\frac {10}3}([b,\infty))}\int_{t_0}^{t_1} \|J(s)|v|^4 v(s)\|_{\frac {10}7} ds
\\
&\le& C|I|^{-\frac 3{10}}|t_1-t_0|^{-\frac 3{10}} \|v\|^4_{L^{10}([t_0,t_1]; L^{10})} \|J(\cdot) v\|_{L^{\frac {10}3}
([t_0,t_1];L^{\frac {10}3})} \\
&\le & C\eta_2^{\frac {10}3} .
\end{eqnarray*}
For the second part, we use Strichartz estimate to control it by
\begin{equation}
C\|J(\cdot)|u|^4 v\|_{L^{\frac {10}7}([t_1,\infty);L^{\frac {10}7})}\le C\|v\|^4_{L^{10}([t_1,\infty);L^{10})}
\|J(\cdot) v\|_{L^{\frac {10}3}([t_1,\infty);L^{\frac {10}3})}.
\label{74}
\end{equation}
At this moment, we follow the same way in proving (\ref{62.1}) to get
$$
\|v\|_{L^{10}([t_1,\infty);L^{10})}\le c(\eta_2),
$$
thus finally, $(\ref{74}) \le c(\eta_2)$. Hence we get Lemma8.1.

Now we are at the position to solve the perturbation problem (\ref{70}). By induction assumption, we see that there exists
constant $C=C(\Lambda_1,\Lambda_2)$ such that
\begin{eqnarray*}
&&\|W\|_{L^{10}([b,\eta_1^4];L^{10})}\le C,\\
&&\|A(\cdot)W\|_{L^{\frac {10}3}([b,\eta_1^4];L^{\frac {10}3})}\le C, A\subset \{ J,H\}.
\end{eqnarray*}
This allows to split $[b,\eta_1^4]$ into finite subintervals
$$
[b,\eta_1^4]=\cup_{j=1}^K I_j=\cup_{j=1}^K [b_{j-1},b_j), b_0=b,
b_K=\eta_1^4,
$$
and such that
$$
\|W\|_{L^{10}(I_j;L^{10})}\sim \nu, \qquad \|A(\cdot)
W\|_{L^{\frac {10} 3}(I_j;L^{\frac {10}3})} \sim \varepsilon.
$$
Then
$$
K\le \max((\frac C{\nu})^{10}, (\frac C{\varepsilon})^{\frac
{10}3}).
$$
If (\ref{70}) has been solved on $[b_0, b_{j-1}]$, and
$$
\|A(b_{j-1})\Gamma(b_{j-1})\|_2\le C^{j-1}c(\eta_2)^{1-\frac
{j-1}{2K}},
$$
then we can solve (\ref{70}) on $[b_{j-1},b_j]$ by proving the solution map
$$
\Phi (\Gamma(t))=U(t-b_{j-1})\Gamma(b_{j-1})-i\int_{b_{j-1}}^t U(t-s)(|v+W+\Gamma|^4(v+W+\Gamma)-|v|^4 v-|W|^4W)(s)ds ,
$$
is contractive on the closed set
\begin{eqnarray*}
&&X:= \biggl \{ \Gamma \in L^{10}(I_j;L^{10}), A(\cdot)\Gamma \in L^{\frac {10}3}(I_j;L^{\frac {10}3}),\mbox{
and}\\
&&\qquad
\|\Gamma\|_X=\|\Gamma\|_{L^{10}(I_j;L^{10})}+\max_{A\in\{J,H\}}\|A(\cdot)\Gamma\|_{L^{\frac
{10}3}(I_j;L^{\frac {10}3})}
 \le C^j c(\eta_2)^{1-\frac j{2K}} \biggr \},
 \end{eqnarray*}
donated with the metric
$$
d(u_1,u_2)=\|u_1-u_2\|_{X},
$$
and complete one step of iteration by estimating
$\|A(b_j)u(b_j)\|_2$ from Strichartz estimate. This is feasible
since we can choose the absolute constants $\varepsilon$, $\nu$,
and the constant $c(\eta_2)$ small enough. The proof is routine
and is omitted. Now we have a finite energy solution $\Gamma$ on
$[b,\eta_1^4]$ such that
\begin{eqnarray*}
\|\Gamma\|^{10}_{L^{10}([b,\eta_1^4];L^{10})}&=& \sum_{j=1}^K \|\Gamma\|^{10}_{L^{10}(I_j;L^{10})}
\\
&\le& \sum_{j=1}^K(C^j c(\eta_2)^{1-\frac j{2K}})^{10} \le C.
\end{eqnarray*}

To conclude the proof of proposition3.5 in case of blow up solution, we collect all the estimates to get
\begin{eqnarray*}
\|u\|_{L^{10}(I^{(3)};L^{10})}&\le&
\|u\|_{L^{10}([b,\eta_1^4];L^{10})}\\
&\le&\|v\|_{L^{10}([b,\eta_1^4];L^{10})}+\|W\|_{L^{10}([b,\eta_1^4];L^{10})}+\|\Gamma\|_{L^{10}([b,\eta_1^4];L^{10})}\\
&\le& C(\Lambda_1 \Lambda_2,\eta_1,\eta_2).
\end{eqnarray*}
Thus, $J_1$ can be controlled by
$$
O\biggl(\frac {C(\Lambda_1
\Lambda_2,\eta_1,\eta_2,\eta_3)}{\eta_1}\biggr)^{10}.
$$
In the same way, $J_2$ also be controlled, and thus
$$
\|u\|_{L^{10}([-\eta_1^4,\eta_1^4];L^{10})}\le C(\Lambda_1,\Lambda_2,\eta_1,\eta_2,\eta_3).
$$
which closes the induction and finally gives proposition3.5.

Finally, we give some comments about this paper. In this paper, we
consider the energy critical Schr\"odinger equation with repulsive
harmonic potential which is quite different from the one without
potential: We have no positive conserved quantity, but decay
quantities which are not time-translation invariant. We solve
these difficulties by first using the decay estimates to reduce
the global problem to a problem on finite time interval, then
completing the analysis by doing induction on a very small
interval.

Now, let's introduce some open problem left by this paper. One
remaining problem is to generalize the result to the higher
dimensional case, which is hopeful in view of the recent work
in\cite{Tao} and will be discussed elsewhere. Another interesting
problem is how to remove the radial assumption. Because in this
case, the equation is not scaling invariance, there is no hope to
follow the same method in \cite{CKSTT}. There are some other
challenging problems concerning the energy critical equation with
focusing nonlinearity and repulsive potential, or defocusing
nonlinearity and attractive potential, which remains completely
open.
\\
\\
\\
\noindent {\bf Acknowledgement:} The author is grateful to
Professor Ping Zhang for introducing this problem and the useful
discussions.

\end{document}